\documentclass[headings=standardclasses]{scrartcl}

\pdfoutput=1

\usepackage[utf8]{inputenc}
\usepackage[T1]{fontenc}
\usepackage[UKenglish]{babel}

\usepackage{graphicx}
\graphicspath{{Plots/}}
\DeclareGraphicsExtensions{.pdf, .jpg, .png, .ps, .eps}
\usepackage{float}

\usepackage{amssymb,amsmath,amsopn,amsthm}
\usepackage{stmaryrd}

\usepackage{hyperref}
\usepackage[backend=bibtex,style=numeric-comp,sorting=none,
  isbn=false,url=false,doi=false,eprint=false]{biblatex}
\usepackage{algpseudocode}
\usepackage{algorithm}
\usepackage{setspace}
\let\Algorithm\algorithm
\renewcommand\algorithm[1][]{\Algorithm[#1]\setstretch{1.2}}
\algdef{SE}[DOWHILE]{Do}{DoWhile}{\algorithmicdo}[1]{\algorithmicwhile\ #1}

\usepackage{pgfplotstable}
\pgfplotsset{width=7cm,compat=newest}
\usepackage{booktabs,colortbl}
\usepackage[font={small}]{caption}
\usepackage{etoolbox}
\makeatletter
\patchcmd{\@maketitle}{\huge}{\LARGE}{}{}
\makeatother

\theoremstyle{plain}
\newtheorem{theorem}{Theorem}

\newtheorem{corollary}{Corollary}

\theoremstyle{definition}
\newtheorem{definition}{Definition}

\newtheorem{remark}{Remark}

\DeclareMathOperator{\vspan}{span}
\DeclareMathOperator{\sons}{sons}
\DeclareMathOperator{\diam}{diam}
\DeclareMathOperator{\dist}{dist}
\DeclareMathOperator*{\argmax}{arg\,max}

\newcommand{\tensor}[1]{\boldsymbol{\mathcal{#1}}}
\newcommand{\norm}[1]{\left\lVert#1\right\rVert}
\newcommand{\seminorm}[1]{\left\lvert#1\right\rvert}
\newcommand{\dotprod}[1]{\left\langle#1\right\rangle}

\begin{document}

\title{Boundary Element Methods for the Wave Equation based on
  Hierarchical Matrices and Adaptive Cross Approximation}

\author{Daniel Seibel
  \thanks{\emph{Department of Mathematics, Saarland University,Germany},
    E-Mail:\href{mailto:seibel@num.uni-sb.de}{\texttt{seibel@num.uni-sb.de}}
  }\footnotemark[1]
}

\date{\today}

\maketitle

\begin{abstract}
  Time-domain Boundary Element Methods (BEM) have been successfully used in
  acoustics, optics and elastodynamics to solve transient problems numerically.
  However, the storage requirements are immense, since the fully populated system
  matrices have to be computed for a large number of time steps or frequencies.
  In this article, we propose a new approximation scheme for the Convolution
  Quadrature Method (CQM) powered BEM, which we apply to scattering problems
  governed by the wave equation. We use $\mathcal{H}^2$-matrix compression in the
  spatial domain and employ an adaptive cross approximation (ACA) algorithm in the
  frequency domain. In this way, the storage and computational costs are reduced
  significantly, while the accuracy of the method is preserved.
\end{abstract}

\section{Introduction}\label{sec:intro}
The numerical solution of wave propagation problems is a crucial task in
computational mathematics and its applications. In this context, BEM play a
special role, since they only require the discretisation of the boundary instead
of the whole domain. Hence, BEM are particularly favourable in situations where
the domain is unbounded, as it is often the case for scattering problems. There,
the incoming wave hits the object and emits a scattered wave, which is to be
approximated in the exterior of the scatterer. 

In contrast to Finite Element or Difference Methods, BEM are based on boundary
integral equations posed in terms of the traces of the solution. For the
classical example of the scalar wave equation, the occurring integral operators
take the form of so called ``retarded potentials'' related to Huygen's principle.
In~\cite{MR859833}, Bamberger and Ha Duong laid the foundation for their analysis
by applying variational techniques in the frequency domain. Since then,
significant improvements have been made, which are explained thoroughly in the
monograph of Sayas~\cite{MR3468871}. Recently, a unified and elegant approach
based on the semi-group theory has been proposed in~\cite{MR3628109}. Besides,
the articles~\cite{ECM2022} and~\cite{MR3628110} give an excellent overview of
the broad topic of time-domain boundary integral equations.

There are three different strategies for the numerical solution of time-dependent
problems with BEM. The classical approach is to treat the time variable
separately and discretise it via a time-stepping scheme. This leads to a sequence
of stationary problems, which can be solved with standard BEM~\cite{MANSUR83}.
However, one serious drawback is the emergence of volume terms even for vanishing
initial conditions and right-hand side. Therefore, additional measures like the
dual-reciprocity method~\cite{MR1144769} are necessary or otherwise the whole
domain needs to be meshed, which undermines the main benefit of BEM.

In comparison to time-stepping methods, space-time methods regard the time
variable as an additional spatial coordinate and discretise the integral
equations directly in the space-time cylinder. To this end, the latter is
partitioned either into a tensor grid or into an unstructured grid made of
tetrahedral finite elements~\cite{LANGER19}. For that reason, space-time methods
feature an inherent flexibility, including adaptive refinement in both time and
space simultaneously as well as the ability to capture moving
geometries~\cite{MR3628108,MR3983080,MR3944570}. However, the computational costs
are high due to the increase in dimensionality and the calculation of the
retarded potentials is far from trivial~\cite{MR4039531}.

Finally, transformation methods like Lubichs' CQM~\cite{MR923707,MR932708}
present an appealing alternative. The key idea is to take advantage of the
convolutional nature of the operators by use of the Fourier-Laplace transform and
to further discretise via linear multi-step~\cite{MR1269502} or Runge-Kutta
methods~\cite{MR2824853,MR2969188}.
Although the transition to the frequency domain comes with certain restrictions,
e.g. the number of time steps has to be fixed a priori, it nevertheless features
some important advantages. Foremost, the approximation involves only spatial
boundary integral operators related to Helmholtz problems. The properties of
these frequency-dependent operators are well studied~\cite{MR1742312} and they
are substantially easier to deal with than retarded potentials. Moreover, the CQM
is applicable for several problems of poro- and visco-elasticity, where only the
Fourier-Laplace transform of fundamental solution is explicitly
known~\cite{MR1489580}. Higher order discretisation spaces~\cite{MR3561482} as
well as variable time step sizes~\cite{MR3520008} are also supported. Apart from
acoustics~\cite{MR3606236,MR3787392,MR3543003}, CQM have been applied
successfully to challenging problems in electrodynamics~\cite{MR3032320},
elastodynamics~\cite{MR2468393,SCHANZ2001,MR3941889} and
quantum mechanics~\cite{MR3628111}.

Regardless of the method in use, we face the same major difficulty: as BEM
typically generate fully populated matrices, the storage and computational costs
are huge. Since this is already valid for the stationary case, so called fast
methods driven by low-rank approximations have been developed for elliptic
equations, see the monographs~\cite{MR2451321,MR2310663,MR3445676}. The crucial
observation is that the kernel function admits a degenerated expansion in the
far-field, which can be exploited by analytic~\cite{MR2121071,MR1954142} or
algebraic compression algorithms~\cite{MR1794343,MR1972724}. In this way, the
numerical costs are lowered to almost linear in the number of degrees of freedom.
The situation becomes even more difficult when moving to time-domain BEM. In the
CQM formulation, several system matrices per frequency need to be assembled,
culminating in a large number of matrices overall. Because they stem
from elliptic problems, it is straightforward to approximate them via standard
techniques~\cite{MR2387910}. Based on the observation that the convolution
weights decay exponentially, cut-off strategies~\cite{MR2470945,MR3292532} have
been developed to accelerate the calculations. Details on how to combine these
two concepts and how to solve the associated systems efficiently are given
in~\cite{MR3267099}. It is also possible to filter out irrelevant frequencies
if a priori information about the solution is known~\cite{MR2452859}.

In this article, we present a novel approach which relies on hierarchical
low-rank approximation in both space and frequency. The main idea is to
reformulate the problem of approximating the convolution weights as a
tensor approximation problem~\cite{pamm.201900072}.
By means of $\mathcal{H}^2$-matrices in space and ACA in frequency, we manage
to reduce the complexity to almost linear in the number of degrees of freedom as
well as in the number of time steps. In other words, the numerical costs are
significantly reduced, which makes the algorithm particularly fast and efficient.

The paper is structured as follows. In Section~\ref{sec:prelim}, we recall the
boundary integral equations and their Galerkin formulation for the wave equation,
which serves as our model problem throughout this article. Subsequently, we
describe the numerical discretisation of the integral equations powered by CQM
and BEM in Section~\ref{sec:numdis}. The next two Sections~\ref{sec:h2matrix}
and~\ref{sec:aca} deal with the low-rank approximation of the associated matrices
and tensors respectively. Afterwards, in Section~\ref{sec:combalgo}, we analyse
the hierarchical approximation and specify the algorithm in its entirety.
Finally, we present numerical examples in Section~\ref{sec:numex} and summarise
our results in Section~\ref{sec:conc}.

\section{Preliminaries}\label{sec:prelim}
\subsection{Formulation of the problem}\label{subsec:formulation}

Let $\Omega^{\textup{in}} \subset \mathbb{R}^3$ be a bounded domain with
Lipschitz boundary $\Gamma = \partial \Omega$ and denote by
$\Omega = \mathbb{R}^3 \setminus \overline{\Omega}_{\textup{in}}$ the exterior
domain. Further, let $n$ be the unit normal vector on $\Gamma$ pointing into
$\Omega$.

\begin{figure}[htb]
  \centering
  \includegraphics{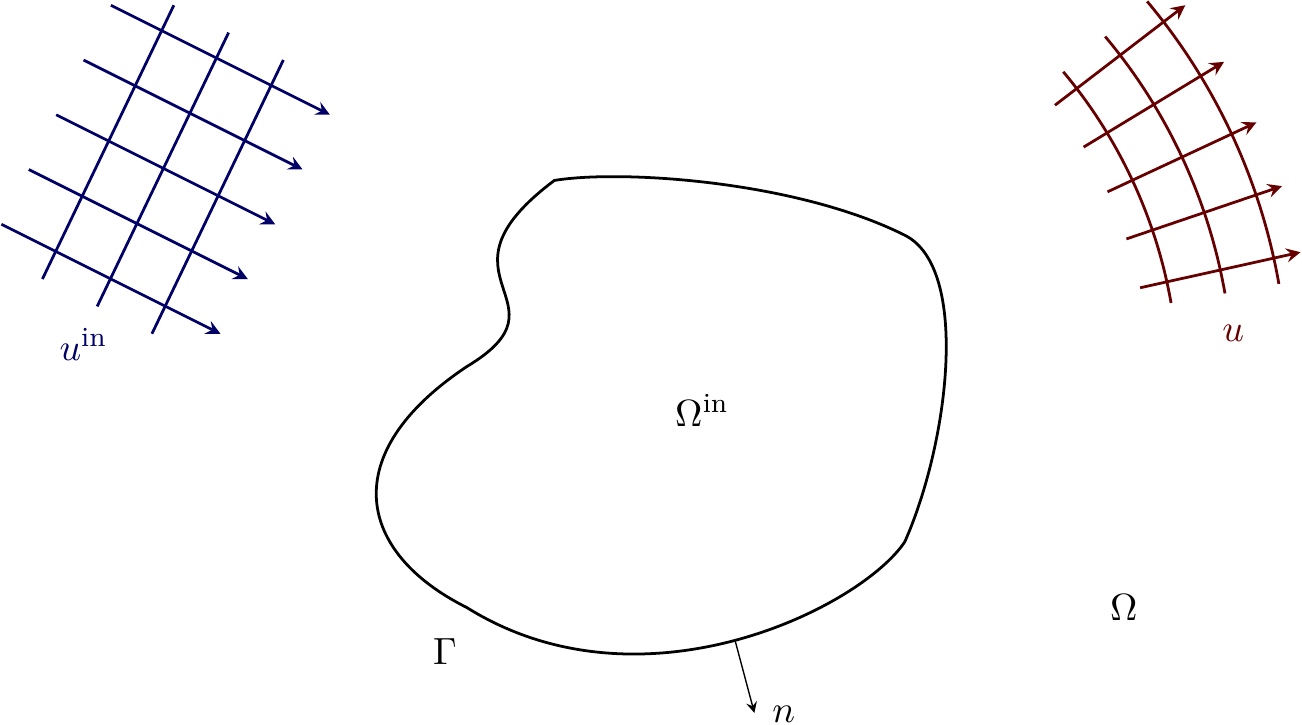}
  \caption{Visualisation of the scattering problem.}~\label{fig:scatter}
\end{figure}

We study the situation depicted in Figure~\ref{fig:scatter} where an
incident wave $u^{\textup{in}}$ is scattered by the stationary obstacle
$\Omega^{\textup{in}}$, causing a scattered wave $u$ to propagate in the free
space $\Omega$. In the absence of external sources, we may assume that $u$
satisfies the homogeneous wave equation
\begin{equation}\label{eq:wave}
  \frac{\partial^2 u}{\partial t^2} (x,t) - c^2 \Delta u (x,t) = 0,
  \quad (x,t) \in \Omega \times (0,T).
\end{equation}
Here, the coefficient $c$ is the speed at which the wave travels in the medium.
Using dimensionless units, we set $c=1$. Moreover, $\Delta$ is the Laplacian in
the spatial domain and $T > 0$ is a fixed final time. Depending on the
characteristics of the scatterer $\Omega^\textup{in}$ and the incoming wave
$u^\textup{in}$, the scattered wave $u$ is subject to boundary conditions posed
on the surface $\Gamma$. We prescribe mixed boundary conditions of the form
\begin{equation}\label{eq:bc}
  \begin{aligned}
    u (x,t) &= g_D (x,t), && (x,t) \in \Gamma_D \times (0,T),\\
    \frac{\partial  u}{\partial n} (x,t) &= g_N (x,t),
    && (x,t) \in \Gamma_N \times (0,T),
  \end{aligned}
\end{equation}
given either on the Dirichlet boundary $\Gamma_D$ or the Neumann boundary
$\Gamma_N$ with
\begin{equation*}
  \Gamma = \overline{\Gamma}_D \cup \overline{\Gamma}_N,
  \quad \Gamma_D \cap \Gamma_N = \varnothing.
\end{equation*}
Furthermore, we assume that $u^{\textup{in}}$ has not reached $\Omega$ yet, which
implies vanishing initial conditions for $u$,
\begin{equation*}
  u (x, 0) = \frac{\partial u}{\partial t} (x,0) = 0, \quad x \in \Omega.
\end{equation*}

\begin{remark}
  Surprisingly enough, regularity results for hyperbolic problems of
  type~\eqref{eq:wave} with non-homogeneous boundary conditions were not
  available until the fundamental works of Lions and Magenes~\cite{MR0350177-9}.
  By the use of sophisticated tools from functional analysis and the theory of
  pseudo-differential operators, they paved the way for the mathematical analysis
  of general second-order hyperbolic systems. Subsequently, their findings were
  substantially improved and we present two examples here. Let
  $\Sigma = \Gamma \times [0,T]$ be the lateral boundary of the space-time
  cylinder $Q = \Omega \times [0,T]$. For the pure Dirichlet case,
  $\Gamma_D = \Gamma$, it is shown in~\cite{MR867669} that
  \begin{equation*}
    g_D \in H^1(\Sigma) \implies  u \in C^0 (H^1(\Omega);[0,T]),\
    \frac{\partial}{\partial n} u \in L^2(\Sigma).
  \end{equation*}
  In comparison, optimal regularity results for the Neumann problem,
  $\Gamma = \Gamma_N$, are derived in~\cite{MR1633000} and are of the form
  \begin{equation*}
    g_N \in L^2(\Sigma) \implies u \in H^{3/4}(Q), \
    u_{| \Sigma} \in H^{1/2}(\Sigma).
  \end{equation*}
  In~\cite{MR1108480} it is shown that this result cannot be improved, i.e.
  that for all $\varepsilon >0$ there exist $g_N \in L^2(\Sigma)$ such that
  $u \not\in H^{3/4 + \varepsilon}(Q)$. From these findings, it becomes evident
  that the situation of mixed conditions like~\eqref{eq:bc} is far from trivial
  and needs special treatment. An alternative approach involves the theory of
  boundary integral equations, which are introduced in Section~\ref{subsec:bie}.
  In combination with the semi-group theory~\cite{MR3628109} or Laplace domain
  techniques~\cite{MR859833,MR3468871}, general transmission problems can be
  treated in a uniform manner. To keep things simple, we refrain from specifying
  the function spaces in the following sections and refer to the given
  publications instead.
\end{remark}

\subsection{Boundary Integral Equations}\label{subsec:bie}
The fundamental solution of~\eqref{eq:wave} is given by
\begin{equation*}
  u^\ast(y - x, t - \tau) =
  \frac{\delta \left(t  - \tau - \seminorm{y - x} \right)}
  {4\pi \seminorm{y - x}}, \quad
  (x,t) \in \Omega \times (0,T),
\end{equation*}
where $\delta$ is the Dirac delta distribution defined by
\begin{equation*}
  \left(\delta(t  - \tau - \seminorm{y - x})\right)(\varphi) =
  \varphi(t - \seminorm{y-x})
\end{equation*}
for smooth test functions $\varphi$. Thus, the behaviour of the
wave at position $x$ and time $t$  is completely determined by its values at
locations $y$ and earlier times $\tau = t - \seminorm{y - x}$. In other words,
an event at $(x,t)$ is only affected by actions that took place on the backward
light cone
\begin{equation*}
  \left\{ (y,\tau) \in \overline{Q} : \tau = t - \seminorm{y-x}\right\}
\end{equation*}
of $(x,t)$ in space-time. Therefore, $u^\ast$ is also known as
\emph{retarded} Green's function and $\tau =t - \seminorm{y - x}$ is called
\emph{retarded} time~\cite{MR0436782}. This property becomes particularly
important in the representation formula,
\begin{equation}\label{eq:representation}
  \begin{aligned}
    u (x,t) &= -\int\limits_0^t \int\limits_\Gamma u^\ast (y - x, t - \tau)\,
    \frac{\partial}{\partial n} u (y, \tau) \mathop{dS(y)}\mathop{d\tau} \\
    &\hspace{12ex} + \int\limits_0^t \int\limits_\Gamma
    \frac{\partial u^\ast}{\partial n_y}(y - x, t - \tau)\, u(y, \tau)
    \mathop{dS(y)}\mathop{d\tau},
  \end{aligned}
\end{equation}
which expresses the solution by the convolution of the boundary data with the
fundamental solution. Although formally specified on the lateral surface
$\Sigma$, the evaluation of the Dirac delta reduces the domain of integration
to its intersection with the backward light cone,
\begin{equation}\label{eq:repret}
  \begin{aligned}
    u (x,t) &= -\int\limits_\Gamma \frac1{4\pi \seminorm{y - x}} \frac{\partial}
    {\partial n} u(y, t - \seminorm{y - x}) \mathop{dS(y)}\mathop{d\tau} \\
    &\hspace{6ex} + \int\limits_\Gamma \frac{(y-x) \cdot n_y}{4 \pi \seminorm{y-x}^2}
    \left( \frac{u(y, t - \seminorm{y - x})}{\seminorm{y - x}} +
      \frac{\partial}{\partial t} u(y, t - \seminorm{y - x}) \right)
    \mathop{dS(y)}\mathop{d\tau}.
  \end{aligned}
\end{equation}
Since the boundary data is only given on $\Gamma_D$ or $\Gamma_N$, we derive a
system of boundary integral equations for its unknown parts.
To this end, we define the trace operators $\gamma_0$ and $\gamma_1$ by
\begin{equation*}
  \gamma_{0,x} v(x) = \lim_{y \to x} v(y), \quad
  \gamma_{1,x} v(x) = \lim_{y \to x} \nabla v(y) \cdot n(x), \quad
  x \in \Gamma, \ y \in \Omega,
\end{equation*}
for sufficiently smooth $v$.
Note that $\gamma_1$ coincides with the normal derivative $\partial/ \partial n$.
Now, we take the Dirichlet trace in~\eqref{eq:representation} and obtain the
first equation
\begin{equation*}
  \begin{aligned}
    \gamma_{0,x} u (x,t) &= - \int\limits_0^t \int\limits_\Gamma u^\ast (y - x, t - \tau)\,
    \gamma_{1,y} u (y, \tau) \mathop{dS(y)}\mathop{d\tau}
    + \frac12 \gamma_{0,x} u (x,t) 
    \\ &\quad - \int\limits_0^t \int\limits_\Gamma \gamma_{1,y} u^\ast (y - x, t - \tau)\,
    \gamma_{0,y} u (y, \tau) \mathop{dS(y)}\mathop{d\tau}
  \end{aligned}
\end{equation*}
for almost every $(x,t) \in \Sigma$. Similarly, the application of the
Neumann trace yields
\begin{equation*}
  \begin{aligned}
    \gamma_{1,x} u (x,t) &= \frac12 \gamma_{1,x} u (x,t) - \int\limits_0^t
    \int\limits_\Gamma \gamma_{1,x} u^\ast (y - x, t - \tau)\, \gamma_{1,x} u (y,\tau)
    \mathop{dS(y)}\mathop{d\tau}\\
    &\quad - \int\limits_0^t \int\limits_\Gamma \gamma_{1,x} \gamma_{1,y}\,
    u^\ast (y - x, t - \tau)\, u (y, \tau) \mathop{dS(y)}\mathop{d\tau}.
  \end{aligned}
\end{equation*}
We identify the terms above with so-called boundary integral operators,
\begin{equation}\label{eq:wavebio}
  \begin{aligned}
    \left( \mathcal{V} w \right) (x,t)
    &= \int\limits_0^t \int\limits_\Gamma u^\ast (y - x, t - \tau)\, w (y, \tau)
    \mathop{dS(y)}\mathop{d\tau},\\
    \left( \mathcal{K} w \right) (x,t)
    &= \int\limits_0^t \int\limits_\Gamma \gamma_{1,y} u^\ast (y - x, t - \tau)\, w (y, \tau)
    \mathop{dS(y)}\mathop{d\tau},\\
    \left( \mathcal{K}^\prime w \right) (x,t)
    &= \int\limits_0^t \int\limits_\Gamma \gamma_{1,x} u^\ast (y - x, t - \tau)\, w (y, \tau)
    \mathop{dS(y)}\mathop{d\tau},\\
    \left( \mathcal{D} w \right) (x,t)
    &= \int\limits_0^t \int\limits_\Gamma \gamma_{1,x} \gamma_{1,y} u^\ast (y - x, t - \tau)
    \, w (y, \tau) \mathop{dS(y)}\mathop{d\tau},\\
    \left( \mathcal{I} w \right) (x,t)
    &= \int\limits_0^t \int\limits_\Gamma \delta(t - \tau - \seminorm{y -x})\, w (y, \tau)
    \mathop{dS(y)}\mathop{d\tau},
  \end{aligned}
\end{equation}
and rewrite the boundary integral equations in matrix form,
\begin{equation}\label{eq:wavebie}
  \begin{pmatrix}
    \gamma_{0,x} u \\ \gamma_{1,x} u
  \end{pmatrix}
  (x,t)
  = \left[
  \begin{pmatrix}
    \tfrac12 \mathcal{I} + \mathcal{K} & -\mathcal{V} \\
    -\mathcal{D} & \tfrac12 \mathcal{I} - \mathcal{K}^\prime
  \end{pmatrix}
  \begin{pmatrix}
    \gamma_{0,y} u \\ \gamma_{1,y} u
  \end{pmatrix}
  \right] (x,t), \quad (x,t) \in \Sigma.
\end{equation}
Hence, the solution to the wave equation~\eqref{eq:wave} is found by solving for
the unknown boundary data in the system above and inserting it into the
representation formula~\eqref{eq:representation}. For the mapping properties of
the operators and the solvability of the equations, we refer to~\cite{MR3628109}.

\subsection{Galerkin formulation}\label{subsec:galform}
We derive a Galerkin formulation of~\eqref{eq:wavebie}. Taking mixed boundary
conditions~\eqref{eq:bc} into account, we choose extensions $\tilde{g}_D$ and
$\tilde{g}_N$ on $\Gamma$ satisfying
\begin{equation*}
  \tilde{g}_D(x,t) = g_D(x,t), \ x \in \Gamma_D, \quad
  \tilde{g}_N(x,t) = g_N(x,t), \ x \in \Gamma_N, \quad t \in [0,T].
\end{equation*}
Furthermore, we decompose the Dirichlet and Neumann traces as follows
\begin{equation*}
  \begin{aligned}
    \gamma_0 u &= \tilde{u} + \tilde{g}_D
    && \textup{with } && \tilde{u} (x,t) = 0 &&
    \textup{ for } (x,t) \in \Gamma_D \times [0,T],\\
    \gamma_1 u &= \tilde{q} + \tilde{g}_N
    && \textup{with } && \tilde{q} (x,t) = 0 &&
    \textup{ for } (x,t) \in \Gamma_N \times [0,T],\\
  \end{aligned}
\end{equation*}
and obtain for $t \in [0,T]$
\begin{equation*}
  \begin{aligned}
    \left( \mathcal{V} \tilde{q} \right) (x,t) -
    \left( \mathcal{K}  \tilde{u} \right) (x,t) &=
    \left[ \Bigl(-\tfrac12 \mathcal{I} + \mathcal{K} \Bigr)
       \tilde{g}_D \right] (x,t) -
    \left( \mathcal{V}  \tilde{g}_N \right) (x,t),
    && x \in \Gamma_D,\\
    \left( \mathcal{K}^\prime  \tilde{q} \right) (x,t) +
    \left( \mathcal{D}  \tilde{u} \right) (x,t) &=
    \left[ \Bigl(-\tfrac12 \mathcal{I} - \mathcal{K}^\prime \Bigr)
       \tilde{g}_N \right] (x,t) -
    \left( \mathcal{D}  \tilde{g}_D \right) (x,t),
    && x \in \Gamma_N.
  \end{aligned}
\end{equation*}
Here, the right hand side is known and we have to solve for the unknown Neumann
data $\tilde{q}$ on the Dirichlet boundary $\Gamma_D$ and the Dirichlet data
$\tilde{u}$ on the Neumann part $\Gamma_N$, respectively. Hence, the Galerkin
formulation is to find $\tilde{q}$ and $\tilde{u}$ such that 
\begin{equation}\label{eq:wavegalform}
  \begin{aligned}
    \dotprod{w,\mathcal{V}\tilde{q}}_{\Gamma_D} -
    \dotprod{w,\mathcal{K}\tilde{u}}_{\Gamma_D} &=
    \dotprod{w,\left( -\tfrac12 \mathcal{I} + \mathcal{K} \right) 
      \tilde{g}_D}_{\Gamma_D} -
    \dotprod{w,\mathcal{V}\tilde{g}_N}_{\Gamma_D},\\
    \dotprod{v,\mathcal{K}^\prime\tilde{q}}_{\Gamma_N} +
    \dotprod{v,\mathcal{D} \tilde{u}}_{\Gamma_N} &=
    \dotprod{v,\left( -\tfrac12 \mathcal{I} - \mathcal{K}^\prime \right) 
      \tilde{g}_N}_{\Gamma_N} -
    \dotprod{v,\mathcal{D}\tilde{g}_D}_{\Gamma_N}
  \end{aligned}
\end{equation}
holds at every time $t \in [0,T]$ and for all test functions $w$ and $v$. The
index of the duality product $\dotprod{\cdot,\cdot}$ indicates on which part
of the boundary it is formed.

\section{Numerical Discretisation}\label{sec:numdis}
In view of the numerical treatment of~\eqref{eq:wavegalform}, we have to
discretise both time and space. We start with the time discretisation.
\subsection{Convolution Quadrature Method}\label{subsec:convquad}
The application of the integral operators in~\eqref{eq:wavebio} requires the
evaluation of the convolution in time of the form
\begin{equation}\label{eq:convolution}
  h(t) = \int\limits_0^t f(t - \tau) g(\tau) \mathop{d\tau},
  \quad 0 \leq t \leq T,
\end{equation}
where $f$ is a distribution and $g$ is smooth. In order to compute such
convolutions numerically, we employ the Convolution Quadrature Method (CQM)
introduced in~\cite{MR923707}. It is based on the Fourier-Laplace transform
defined by
\begin{equation*}
  \hat{f}(s) = \int\limits_0^\infty f(t) e^{-st} \mathop{dt} \in \mathbb{C},
  \quad s \in \mathbb{C},
\end{equation*}
and on the observation that we can replace $f$ by the inverse transform of
$\hat{f}$ and change the order of integration, i.e.
\begin{equation*}
  h(t) = \int\limits_\mathcal{C} \int\limits_0^t e^{(t-\tau)s}
  g(\tau) \mathop{d\tau}
  \hat{f}(s) \mathop{ds},
  \quad 0 \leq t \leq T.
\end{equation*}
The integration is performed along the contour
\begin{equation*}
  \mathcal{C} = \left\{ s \in \mathbb{C} : \sigma + \imath z,\
  z\in\mathbb{R}\right\},
\end{equation*}
where $\sigma >0$ is greater than the real part of all singularities of
$\hat{f}$. In further steps, the inner integral is approximated by a linear
multi-step method and Cauchy's integral formula is used. In this way, the CQM
yields approximations of~\eqref{eq:convolution} at discrete time points
\begin{equation*}
  t_n = n \Delta t,\quad \Delta t = T / N, \quad n = 0,\ldots,N,
\end{equation*}
via the quadrature formula
\begin{equation}\label{eq:cqm}
  h (t_n) \approx \sum_{k=0}^n \omega_{n-k}\, g(t_k),
  \quad n=0,\ldots,N.
\end{equation}
For a parameter $R>0$, the quadrature weights are defined by
\begin{equation}\label{eq:omega_n}
  \omega_n = \frac{R^{-n}}{N} \sum_{\ell = 0}^{N-1} \hat{f}
  \left( s_\ell \right) e^{\tfrac{-2\pi \imath}{N} n\ell}, \quad
  n = 0,\ldots, N,
\end{equation}
and they are given in terms of $\hat{f}$ sampled at specific frequencies
\begin{equation}\label{eq:cqm_freq}
  s_\ell = \frac{\chi
    \left(R\cdot e^{\tfrac{2\pi \imath}{N} \ell}\right)}{\Delta t} \in \mathbb{C},
  \quad \ell = 0, \ldots, N-1,
\end{equation}
which depend on the characteristic function $\chi$ of the multi-step
method~\cite{MR1227985}. For the choice of the parameter $R$ and further
details we refer the reader to~\cite{MR923707,MR932708,MR1269502}.

Returning to the setting of~\eqref{eq:wavegalform}, we apply the CQM to
approximate the expressions occurring in the Galerkin formulation, for instance
$h(t) = \dotprod{w,\mathcal{V}\tilde{q}(t)}_{\Gamma_D}$, at equidistant time
steps $t = t_n$. For the the single layer operator $\mathcal{V}$ we
obtain
\begin{equation*}
  \begin{aligned}
    \dotprod{w,\mathcal{V}\tilde{q}(t_n)}_{\Gamma_D}
    &= \int\limits_0^{t_n} \int\limits_{\Gamma_D}
    \int\limits_{\Gamma} u^\ast (y-x, t_n - \tau)\,\tilde{q}(y,\tau)
    \mathop{dS(y)} w(x) \mathop{dS(x)} \mathop{d\tau}\\
    &\approx \int\limits_{\Gamma_D} \int\limits_{\Gamma} \sum_{k=0}^n
    \omega_{n-k}(y-x)\, \tilde{q}(y,t_k) \mathop{dS(y)} w(x)
    \mathop{dS(x)},
  \end{aligned}
\end{equation*}
with quadrature weights
\begin{equation*}
  \omega_{n-k}(y-x) = \frac{R^{-(n-k)}}{N} \sum_{\ell = 0}^{N-1} \hat{u}^\ast
  \left(y-x, s_\ell\right) e^{\tfrac{-2\pi \imath}{N} (n-k)\ell}.
\end{equation*}
The transformed fundamental solution $\hat{u}^\ast$ is precisely the fundamental
solution of the Helmholtz equation for complex frequencies $s \in \mathbb{C}$,
\begin{equation*}
  -\Delta \hat{u} + s^2 \hat{u} = 0,
\end{equation*}
and has the representation
\begin{equation}\label{eq:fundsol}
  \hat{u}^\ast (y-x, s) = \frac{e^{\textstyle{-s \seminorm{y-x}}}}
  {4 \pi \seminorm{y-x}}.
\end{equation}
In contrast to the retarded fundamental solution $u^\ast$, its transform
$\hat{u}^\ast$ defines a smooth function for $x \neq y$ and all $s$. Hence, the
CQM formulation has the advantage that distributional kernel functions are
avoided. We set $\tilde{q}_k = \tilde{q}(\cdot,t_k)$ and introduce the operator
\begin{equation}\label{eq:vhat}
  \widehat{\mathbf{V}}_{n-k}\, \tilde{q}_k (x) =
  \frac{R^{-(n-k)}}{N} \sum_{\ell = 0}^{N-1} e^{\tfrac{-2\pi \imath}{N} (n-k)\ell}
  \int\limits_{\Gamma} \hat{u}^\ast \left(y-x, s_\ell \right)
  \tilde{q}_k(y) \mathop{dS(y)},
\end{equation}
which acts only on the spatial component. Finally, we end up with the
approximation
\begin{equation*}
  \dotprod{w, \mathcal{V}\tilde{q}(t_n)}_{\Gamma_D} \approx
  \sum_{k=0}^n \dotprod{w, \widehat{\mathbf{V}}_{n-k}\,\tilde{q}_k}_{\Gamma_D},
\end{equation*}
where the continuous convolution is now replaced by a discrete one.
Repeating this procedure for the other integral operators leads to the
time-discretised Galerkin formulation: find $\tilde{u}_k$ and $\tilde{q_k}$,
$k=0,\ldots,N$, such that
\begin{equation}\label{eq:cqmgalform}
  \begin{aligned}
    \sum_{k=0}^n \left(
      \begin{aligned}
        \dotprod{w,\widehat{\mathbf{V}}_{n-k}\, \tilde{q}_k}_{\Gamma_D} -
        \dotprod{w,\widehat{\mathbf{K}}_{n-k}\, \tilde{u}_k}_{\Gamma_D} \\
        \dotprod{v,\widehat{\mathbf{K}}^\prime_{n-k}\, \tilde{q}_k}_{\Gamma_N} +
        \dotprod{v,\widehat{\mathbf{D}}_{n-k}\, \tilde{u}_k}_{\Gamma_N}
      \end{aligned}
    \right) &=\\
    &\hspace{-10em}\sum_{k=0}^n \left(
      \begin{aligned}
        \dotprod{w,\left( -\tfrac12 \widehat{\mathbf{I}}_{n-k} +
            \widehat{\mathbf{K}}_{n-k} \right) \tilde{g}_{D,k}}_{\Gamma_D} -
        \dotprod{w,\widehat{\mathbf{V}}_{n-k}\, \tilde{g}_{N,k}}_{\Gamma_D} \\
        \dotprod{v,\left( -\tfrac12 \widehat{\mathbf{I}}_{n-k} -
            \widehat{\mathbf{K}}^\prime_{n-k}\right)\tilde{g}_{N,k}}_{\Gamma_N}-
        \dotprod{v,\widehat{\mathbf{D}}_{n-k} \tilde{g}_{D,k}}_{\Gamma_N}
      \end{aligned}
    \right)
  \end{aligned}
\end{equation}
holds for all test functions $w$ and $v$ and $n=0, \ldots, N$. 

From its Definition~\eqref{eq:vhat}, we see that the single layer operator
$\widehat{\mathbf{V}}_{n-k}$ admits the representation
\begin{equation}\label{eq:slphelm}
  \widehat{\mathbf{V}}_{n-k}\, \tilde{q}_k (x) =
  \frac{R^{-(n-k)}}{N} \sum_{\ell = 0}^{N-1}
  e^{\tfrac{-2\pi \imath}{N} (n-k)\ell}
  \left(\mathbf{V}_\ell\, \tilde{q}_k \right) (x)
\end{equation}
as a scaled discrete Fourier transform of operators $\mathbf{V}_\ell$,
\begin{equation*}
  \mathbf{V}_\ell\, \tilde{q}_k (x) = \int\limits_{\Gamma}
  \hat{u}^\ast \left(y-x, s_\ell \right) \tilde{q}_k(y) \mathop{dS(y)},
  \quad \ell = 0, \ldots, N-1.
\end{equation*}
These are exactly the single layer operators corresponding to the Helmholtz
equations with frequencies $s_\ell$. In the same manner, the other integral
operators $\widehat{\mathbf{K}}_{n-k}$, $\widehat{\mathbf{K}}^\prime_{n-k}$ and
$\widehat{\mathbf{D}}_{n-k}$ may be written in terms of the respective operators
$\mathbf{K}_\ell$, $\mathbf{K}_\ell^\prime$, $\mathbf{D}_\ell$ of the Helmholtz
equations. Since $\textup{Im} (s_\ell) \neq 0$, the single layer operator
$\mathbf{V}_\ell$ and hypersingular operator $\mathbf{D}_\ell$ are
elliptic~\cite{MR1822275,MR859833}. The operator $\widehat{\mathcal{I}}_{n-k}$
on the other hand can be easily calculated,
\begin{equation*}
  \widehat{\mathbf{I}}_{0} = \mathcal{I}, \quad
  \widehat{\mathbf{I}}_{n-k} = 0,\quad k \neq n.
\end{equation*}
Therefore, the spatial discretisation of~\eqref{eq:cqmgalform} is equivalent to
a spatial discretisation of a sequence of Helmholtz problems.

\subsection{Galerkin Approximation}\label{subsec:galapx}
We assume that the boundary $\Gamma$ admits a decomposition into
flat triangular elements, which belong either to $\Gamma_D$ or $\Gamma_N$.
We define boundary element spaces of constant and linear order
\begin{equation*}
  S_h^0(\Gamma_D) = \vspan {\left\{ \varphi_m^0 \right\}}_{m=1}^{M_D}, \quad
  S_h^1(\Gamma_N) = \vspan {\left\{ \varphi_m^1 \right\}}_{m=1}^{M_N},
\end{equation*}
and global variants
\begin{equation*}
  S_h^0(\Gamma) = \vspan {\left\{ \varphi_m^0 \right\}}_{m=1}^{M_0}, \quad
  S_h^1(\Gamma) = \vspan {\left\{ \varphi_m^1 \right\}}_{m=1}^{M_1}.
\end{equation*}
Then, we follow the ansatz
\begin{equation*}
  \tilde{q}_n = \sum_{m=1}^{M_D} \underline{q}_n[m]\, \varphi_m^0
  \in S_h^0(\Gamma_D), \quad
  \tilde{u}_n = \sum_{m=1}^{M_N} \underline{u}_n[m]\, \varphi_m^1
  \in S_h^1(\Gamma_N),\quad
  n=0,\ldots,N,
\end{equation*}
with $\underline{q}_n\in\mathbb{R}^{M_D}$, $\underline{u}_n\in\mathbb{R}^{M_N}$
to approximate the unknown Neumann and Dirichlet data.
Likewise, the boundary conditions are represented by coefficient vectors
$\underline{g}_n^N\in\mathbb{R}^{M_0}$ and $\underline{g}_n^D\in\mathbb{R}^{M_1}$,
which are determined by $L^2$-projections onto $S_h^i(\Gamma)$.

As pointed out before, we begin with the discretisation of the boundary integral
operators of Helmholtz problems. For each frequency $s_\ell$, $\ell=0,\ldots,N-1$,
we have boundary element matrices
\begin{equation*}
  \begin{gathered}
    V_\ell \in \mathbb{C}^{M_D \times M_0},\quad
    K_\ell \in \mathbb{C}^{M_D \times M_1},\quad
    K^\prime_\ell \in \mathbb{C}^{M_N \times M_0},\quad
    D_\ell \in \mathbb{C}^{M_N \times M_1},\\
    I \in \mathbb{R}^{M_D \times M_1}, \quad
    I^\prime \in \mathbb{R}^{M_N \times M_0},
  \end{gathered}
\end{equation*}
defined by
\begin{equation*}
  \begin{aligned}
    V_\ell[m,i] &= \dotprod{\varphi_m^0,\mathbf{V}_\ell\varphi_i^0}_{\Gamma_D},&
    K_\ell[m,j] &= \dotprod{\varphi_m^0,\mathbf{K}_\ell\varphi_j^1}_{\Gamma_D},\\
    K^\prime_\ell[p,i] &= \dotprod{\varphi_p^1,\mathbf{K}^\prime_\ell\varphi_i^0
    }_{\Gamma_N},&
    D_\ell[p,j] &= \dotprod{\varphi_p^1,\mathbf{D}_\ell\varphi_j^1}_{\Gamma_N},\\
    I[m,j] &= \dotprod{\varphi_m^0,\varphi_j^1}_{\Gamma_D},&
    I^\prime[p,i] &= \dotprod{\varphi_p^1,\varphi_i^0}_{\Gamma_N},
  \end{aligned}
\end{equation*}
with $i=1,\ldots,M_0$, $j=1,\ldots,M_1$, $m=1,\ldots,M_D$, $p=1,\ldots,M_N$.
Just as in~\eqref{eq:slphelm}, these auxiliary matrices are then transformed
to obtain the integration weights of the CQM. In the case of the single layer
operators this amounts to
\begin{equation}\label{eq:dft}
  \widehat{V}_n = \frac{R^{-n}}{N}
  \sum_{\ell = 0}^{N-1} e^{\tfrac{-2\pi \imath}{L} n\ell}\ V_\ell,\quad
  n = 0, \ldots, N,
\end{equation}
such that
\begin{equation*}
  \widehat{V}_n[m,i] = \dotprod{\varphi_m^0,\widehat{\mathbf{V}}_n \varphi_i^0
  }_{\Gamma_D}, \quad i=1,\ldots,M_0,\quad m=1,\ldots,M_D,
\end{equation*}
holds by linearity.
Moreover, we identify sub-matrices
\begin{equation*}
  \begin{aligned}
    \widehat{V}^D_{n}&=\widehat{V}_n[1:M_D,1:M_D]
    \in \mathbb{C}^{M_D \times M_D}, \\
    \widehat{K}^N_{n}&=\widehat{K}_n[1:M_D,1:M_N]
    \in \mathbb{C}^{M_D \times M_N}, \\
    \widehat{D}^N_{n}&=\widehat{D}_n[1:M_N,1:M_N]
    \in \mathbb{C}^{M_N \times M_N}.
  \end{aligned}
\end{equation*}
The Galerkin approximation of~\eqref{eq:cqmgalform} is then equivalent to the
system of linear equations
\begin{equation}\label{eq:cqmgalapx}
  \begin{pmatrix}
    \widehat{V}^D_0 & -\widehat{K}^N_0 \\[1ex]
    (\widehat{K}^N_0)^\top & \widehat{D}^N_0
  \end{pmatrix}
  \begin{pmatrix}
    \underline{q}_n \\[1ex] \underline{u}_n
  \end{pmatrix}
  =
  \begin{pmatrix}
    f_n^D \\[1ex] f_n^N
  \end{pmatrix},
  \quad n=0,\ldots,N,
\end{equation}
with right-hand side
\begin{equation}\label{eq:cqmrhs}
  \begin{aligned}
    \begin{pmatrix}
      f_n^D \\[1ex] f_n^N
    \end{pmatrix}
    &=-\frac12
    \begin{pmatrix}
      I \underline{g}_n^D \\[1ex]  I^\prime \underline{g}_n^N
    \end{pmatrix}
    +
    \sum_{k=0}^n
    \begin{pmatrix}
      -\widehat{V}_{n-k} & \widehat{K}_{n-k} \\[1ex]
      -\widehat{K}^\prime_{n-k} & -\widehat{D}_{n-k}
    \end{pmatrix}
    \begin{pmatrix}
      \underline{g}_k^N \\[1ex] \underline{g}_k^D
    \end{pmatrix}\\
    &\hspace{20ex} +
    \sum_{k=0}^{n-1}
    \begin{pmatrix}
      -\widehat{V}^D_{n-k} & \widehat{K}^N_{n-k} \\[1ex]
      -(\widehat{K}^N_{n-k})^\top & -\widehat{D}^N_{n-k}
    \end{pmatrix}
    \begin{pmatrix}
      \underline{q}_{k} \\[1ex] \underline{u}_{k}
    \end{pmatrix}.
  \end{aligned}
\end{equation}
While the first row in~\eqref{eq:cqmrhs} corresponds to the right-hand side
of~\eqref{eq:cqmgalform}, the second row contains the boundary values of the
previous time steps and results from the convolutional structure of the CQM
approximation. Since the left-hand side of the linear system stays the same for
every time step, only one matrix inversion has to be performed throughout the
whole simulation. To be more precise, system~\eqref{eq:cqmgalapx} is equivalent
to the decoupled system
\begin{equation*}
  \begin{aligned}
    S \underline{u}_n &= f_n^N - (\widehat{K}^N_0)^\top
    (\widehat{V}^D_0)^{-1} f_n^D,\\
    \widehat{V}^D_0 \underline{q}_n &= \widehat{K}^N_0 \underline{u}_n +
    f_n^D,
  \end{aligned}
\end{equation*}
where $S$ is the Schur complement
\begin{equation*}
  S = \widehat{D}^N_0 + (\widehat{K}^N_0)^\top (\widehat{V}^D_0)^{-1}
  \widehat{K}^N_0.
\end{equation*}
Since both $\widehat{V}_0$ and $S$ are real symmetric and positive definite,
we factorise them via LU decomposition once for $n=0$ and use forward and
backward substitution to solve the systems progressively in time.

Naturally, the assembly of the boundary element matrices and the computation of
the right-hand side~\eqref{eq:cqmrhs} for each step are the most demanding
parts of the algorithm, both computational and storage wise. Due to the fact that
the matrices are generally fully populated, sparse approximation techniques are
indispensable for large scale problems. Compared to stationary problems, this is
even more crucial here, as the amount of numerical work scales with the number of
time steps. Therefore, it is necessary to not only approximate in the spatial but
also in the frequency variable. It proves to be beneficial to interpret the array
of matrices $V_k$, $k = 0, \ldots, N-1$, as a third order tensor
\begin{equation}\label{eq:tensor}
  \tensor{V}[i,j,k] = V_k[i,j].
\end{equation}
In this way, we can restate the problem within the frame of general tensor
approximation and compression. To that end, we introduce low-rank factorisations
which make use of the tensor product.
\begin{definition}[Tensor Product]\label{def:tenprod}
  For matrices $A^{(j)} \in \mathbb{C}^{r_j \times I_j}$, $j = 1, \ldots, d$ and
  a tensor $\tensor{X} \in \mathbb{C}^{I_1 \times \cdots \times I_d}$,
  we define the tensor or mode product $\times_j$ by
  \begin{equation*}
    \left(\tensor{X} \times_j A^{(j)} \right)
    [i_1,\ldots,i_{j-1},\ell,i_{j+1},\ldots,i_d] =
    \sum_{i_j = 1}^{I_j} \tensor{X}[i_1,\ldots,i_d] \, A^{(j)}[\ell,i_j],
    \quad \ell = 1, \ldots, r_j.
  \end{equation*}
\end{definition}
Because of the singular nature of the fundamental solution, a global low-rank
approximation of $\tensor{V}$ is practically not achievable. Instead, we follow
a hierarchical approach where we partition the tensor into blocks, which we
approximate individually. Our scheme is based on $\mathcal{H}^2$-matrix
approximation in the spatial domain, i.e. in $i$ and $j$, and ACA in the
frequency, i.e. in $k$.

\section{Hierarchical Matrices}\label{sec:h2matrix}
The boundary element matrices $V_k$ are of the form
\begin{equation*}
  G[i,j] = \int\limits_\Gamma \int\limits_\Gamma g(x,y)
  \varphi_j(y) \mathop{dS(y)} \psi_i(x) \mathop{dS(x)}, \quad i \in I, j \in J,
\end{equation*}
where $\psi_i$ and $\varphi_j$ are trial and test functions respectively with
index sets $I$ and $J$ and $g$ is a kernel function.
We associate with each $i\in I$ and $j \in J$ sets $X_i$ and $Y_j$, which
correspond to the support of the trial and test functions $\psi_i$ and
$\varphi_j$. For $r \subset I$ and $c \subset J$, we define
\begin{equation*}
  X_r = \bigcup_{i \in r}\, X_i \textup{ and } Y_c = \bigcup_{j \in c}\, Y_j.
\end{equation*}
Moreover, we choose axis-parallel boxes $B_r$ and $B_c$ that contain the sets
$X_r$ and $Y_c$, respectively.

Since $g$ is non-local, the matrix $G$ is typically fully populated.
However, if $X_r$ and $Y_c$ are well separated, i.e. if they satisfy the
admissibility condition
\begin{equation}\label{eq:admissibility}
  \max \left\{ \diam(X_r), \diam(X_c) \right\} \leq \eta \dist(X_r,X_c)
\end{equation}
for fixed $\eta>0$, then the kernel function admits the degenerated
expansion
\begin{equation}\label{eq:degexp}
  g(x,y) \approx \sum_{\mu=1}^p \sum_{\nu=1}^p L_{r,\mu}(x)
  g(\xi_{r,\mu},\xi_{c,\nu}) L_{c,\nu}(y), \quad
  x \in X_r, y \in Y_c,
\end{equation}
into Lagrange polynomials $L_{r,\mu}$ and $L_{c,\nu}$ on $X_r\times Y_c$ . Here,
we choose tensor products $\xi_{r,\mu}$ and $\xi_{c,\nu}$ of Chebyshev points in
$B_r$ and $B_c$ as interpolation points. In doing so, the double integral reduces
to a product of single integrals
\begin{equation*}
  G[i,j] \approx \sum_{\mu=1}^p \sum_{\nu=1}^p
  g(\xi_{r,\mu}, \xi_{c,\nu})
  \int\limits_\Gamma L_{r,\mu}(x) \psi_i(x) \mathop{dS(x)}
  \int\limits_\Gamma L_{c,\nu}(y) \varphi_j(y) \mathop{dS(y)},
  \ (i,j) \in r \times c,
\end{equation*}
which results in the low-rank approximation of the sub-block
\begin{equation}\label{eq:uniformlr}
  G[b] \approx U_b\, S_b\, W_b^H, \quad b = r \times c,
\end{equation}
with
\begin{equation*}
  \begin{aligned}
    U_b[i,\mu] &= \int\limits_\Gamma L_{r,\mu}(x) \psi_i(x) \mathop{dS(x)},
    && i \in r,\quad \mu = 1, \ldots, p,\\
    S_b[\mu,\nu] &= g(\xi_{r,\mu}, \xi_{c,\nu}),
    && \mu,\nu = 1, \ldots, p,\\
    W_b[j,\nu] &= \int\limits_\Gamma L_{c,\nu}(y) \varphi_j(y) \mathop{dS(y)},
    && j \in c,\quad \nu = 1, \ldots, p.
  \end{aligned}
\end{equation*}
By approximating suitable sub-blocks with low-rank matrices, we obtain a
hierarchical matrix approximation of $G$. This approach leads to a reduction of
both computational and storage costs for assembling $G$ from quadratic to almost
linear in $\#I$ and $\#J$, where $\#I$ denotes the cardinality of the set $I$.

\subsection{Matrix Partitions}\label{subsec:partitions}
In the following, we give a short introduction on hierarchical matrices based
on the monographs~\cite{MR2451321,MR2767920}. Since only sub-blocks that satisfy
the admissibility condition~\eqref{eq:admissibility} permit accurate low-rank
approximations, a partition of the matrix indices $I \times J$ into admissible
and inadmissible blocks is required. To this end, we define cluster trees for
$I$ and $J$.

\begin{definition}[Cluster trees]
  Let $\mathcal{T}(I)$ be a tree with nodes $\varnothing \neq r \subset I$.
  We call $\mathcal{T}(I)$ a cluster tree if the following conditions hold:
  \begin{enumerate}
  \item $I$ is the root of $\mathcal{T}(I)$.
  \item If $r \in \mathcal{T}(I)$ is not a leaf, then $r$ is the disjoint union
    of its sons
    \begin{equation*}
      r = \bigcup_{r^\prime \in \sons (r)} r^\prime.
    \end{equation*}
  \item $\# \sons(r) \neq 1$ for $r \in \mathcal{T}(I)$.
  \end{enumerate}
  We denote by $\mathcal{L}(\mathcal{T}(I))$ the set of leaf clusters
  \begin{equation*}
    \mathcal{L}(\mathcal{T}(I)) = \left\{ r \in \mathcal{T}(I) : \sons(r)
      = \varnothing \right\}.
  \end{equation*}
  Moreover, we assume that the size of the clusters is bounded from below, i.e.
  \begin{equation*}
    \# r > n_{\min} > 1, \quad r \in \mathcal{T}(I),
  \end{equation*}
  in order to control the number of clusters and limit the overhead in practical
  applications.
\end{definition}

There are several strategies to perform the clustering efficiently. For
instance, the geometric clustering in~\cite{MR3445676} constructs the cluster
tree recursively by splitting the bounding box in the direction with largest
extent.
Alternatively, the principal component analysis can be used to produce
well-balanced cluster trees~\cite{MR2451321}.

Since searching the whole index set $I \times J$ for an optimal partition is not
reasonable, we restrict ourselves to partitions which are based on cluster
trees of $I$ and $J$.

\begin{definition}[Block cluster trees]
  Let $\mathcal{T}(I)$ and $\mathcal{T}(J)$ be cluster trees. We construct the
  block cluster tree $\mathcal{T}(I \times J)$ by
  \begin{enumerate}
  \item setting $I \times J$ as the root of $\mathcal{T}(I \times J)$,
  \item and defining the sons recursively starting with $r \times c$ for
    $r = I$ and $c = J$:
    \begin{equation*}
      \sons(r \times c) = \left\{
        \begin{aligned}
          & \sons(r) \times c, && \textup{if} \sons(r) \neq \varnothing
          \textup{ and } \sons(c) = \varnothing,\\
          & r \times \sons(c), && \textup{if} \sons(r) = \varnothing
          \textup{ and } \sons(c) \neq \varnothing,\\
          & \sons(r) \times \sons(c), && \textup{if} \sons(r) \neq \varnothing
          \textup{ and } \sons(c) \neq \varnothing,\\
          & \varnothing, && \textup{if } r \times c \textup{ is admissible or}\\
          &&&\sons(r) = \sons(c) = \varnothing.
        \end{aligned}
      \right.
    \end{equation*}
  \end{enumerate}
\end{definition}

Then, the set of leaves $\mathcal{L}(\mathcal{T}(I \times J))$ is a partition in
the following sense.

\begin{definition}[Admissible partition]
  We call $\mathcal{P}$ a partition of $I \times J$ with respect to the block
  cluster tree $\mathcal{T}(I \times J)$ if
  \begin{enumerate}
  \item $\mathcal{P} \subset \mathcal{T}(I \times J)$,
  \item $b, b^\prime \in \mathcal{P} \implies b \cap b^\prime = \varnothing
    \textup{ or } b = b^\prime$,
  \item $\dot{\bigcup}_{b \in \mathcal{P}}\, b = I \times J$.
  \end{enumerate}
  Moreover, $\mathcal{P}$ is said to be admissible if every
  $r \times c \in \mathcal{P}$ is either admissible~\eqref{eq:admissibility} or
  \begin{equation*}
    \max \left\{\#r, \#c \right\} \leq n_{\min}.
  \end{equation*}
  In this case, the near and far field of $\mathcal{P}$ are defined by
  \begin{equation*}
    \mathcal{P}^- = \left\{ r \times c \in \mathcal{P}
      : \max \left\{\#r, \#c \right\} \leq n_{\min} \right\}, \quad
    \mathcal{P}^+ = \mathcal{P} \setminus \mathcal{P}^-.
  \end{equation*}
\end{definition}

\begin{figure}[htb]
  \centering
  \includegraphics[width=0.5\linewidth]{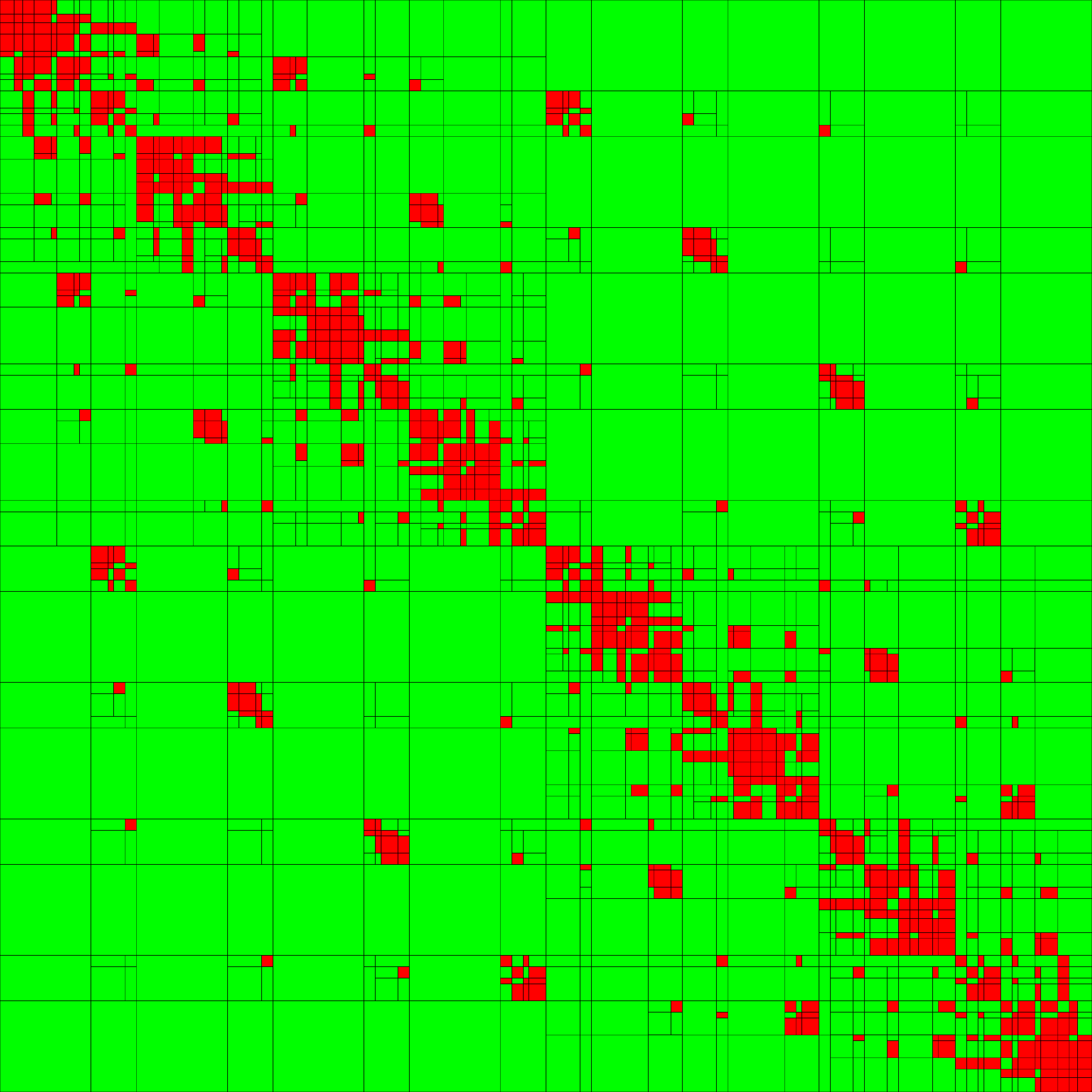}
  \caption{Visualisation of a matrix partition. Green blocks are admissible,
    whereas red ones are inadmissible.}\label{fig:partition}
\end{figure}

Thus, the near field $\mathcal{P}^-$ describes those blocks of $G$ that are
stored in full, because they are inadmissible or simply too small. On the other
hand, the far field $\mathcal{P}^+$ contains admissible blocks only, which are
approximated by low-rank matrices. In Figure~\ref{fig:partition} a partition for
the single layer potential is visualised.

\begin{remark}
  Since evaluating the admissibility condition~\eqref{eq:admissibility} is rather
  expensive, we use the alternative condition
  \begin{equation}\label{eq:qadm}
    \max \left\{ \diam(B_r), \diam(B_c) \right\} \leq \eta \dist(B_r,B_s),
  \end{equation}
  which operates on the bounding boxes and is easier to check.
\end{remark}

\subsection{$\mathcal{H}^2$-matrices}\label{subsec:h2}

One special class of hierarchical matrices consists of $\mathcal{H}^2$-matrices.
They are based on the observation that the matrices $U_b$ and $W_b$ in the
low-rank factorisation~\eqref{eq:uniformlr} of the far field block $b =r\times c$
only depend on the respective row cluster $r$ or column cluster $c$ and not on
the block $b$ itself.

\begin{definition}[$\mathcal{H}^2$-matrices]
  Let $\mathcal{P}$ be an admissible partition of $I \times J$.
  \begin{enumerate}
  \item We call
    \begin{equation*}
      {\left( U_r \right)}_{r \in \mathcal{T}(I)},
      \quad U_r \in \mathbb{C}^{r \times k_r}, \quad k_r > 0,
    \end{equation*}
    (nested) cluster basis, if for all non-leaves $r \in \mathcal{T}(I) \setminus
    \mathcal{L}(\mathcal{T}(I))$ transfer matrices
    \begin{equation*}
      E_{r^\prime,r} \in \mathbb{C}^{k_{r^\prime} \times k_r},
      \quad r^\prime \in \sons(r),\     
    \end{equation*}
    exist such that
    \begin{equation*}
      U_r =
      \begin{pmatrix}
        U_{r_1} E_{r_1,r} \\
        \dots \\
        U_{r_p} E_{r_p,r}
      \end{pmatrix},
      \quad \sons(r) = \left\{ r_1, \ldots, r_p \right\}.
    \end{equation*}
  \item $G$ is called $\mathcal{H}^2$-matrix with row cluster basis
    ${\left( U_r \right)}_{r \in \mathcal{T}(I)}$ and column cluster basis
    ${\left( W_c \right)}_{c \in \mathcal{T}(J)}$, if for each far field block
    $b = r \times c$ there are coupling matrices
    $S_b \in \mathbb{C}^{k_r \times k_c}$ such that
    \begin{equation*}
      G[b] = U_r\, S_b\, W_c^H.
    \end{equation*}    
  \end{enumerate}
\end{definition}

In view of our interpolation scheme, we observe that the Lagrange polynomials
of the father cluster $r \in \mathcal{T}(I)$ can be expressed by the Lagrange
polynomials of its son clusters $r^\prime \in \sons(r)$ via interpolation,
\begin{equation*}
  {L_{r,\mu}}(x) = \sum_{\lambda = 1}^p
  L_{r,\mu}(\xi_{r^\prime,\lambda}) L_{r^\prime,\lambda}(x),
  \quad x \in B_{r^\prime}.
\end{equation*}
Hence, by choosing transfer matrices
\begin{equation*}
  E_{r^\prime,r}[\lambda,\mu] = L_{r,\mu}(\xi_{r^\prime,\lambda}),
\end{equation*}
the cluster basis becomes nested
\begin{equation*}
  \begin{aligned}
    U_r[i,\mu] &= \int\limits_{X_i} L_{r,\mu}(x) \psi_i(x) \mathop{dS(x)}
    = \sum_{\lambda = 1}^p L_{r,\mu}(\xi_{r^\prime,\lambda}) \int\limits_{X_i}
    L_{r^\prime,\lambda}(x) \psi_i(x) \mathop{dS(x)} \\
    &= \sum_{\lambda =1}^p E_{r^\prime,r}[\lambda,\mu]\ U_{r^\prime}[i,\lambda]
    = \left(U_{r^\prime} E_{r^\prime,r} \right)[i,\mu], \quad i \in r^\prime.
  \end{aligned}
\end{equation*}
Algorithm~\ref{algo:h2matrix} describes the assembly of cluster bases and
summarises the construction of an $\mathcal{H}^2$ -matrix by interpolation.

In the following, let $\widetilde{G}$ be the $\mathcal{H}^2$-approximation of the
dense Galerkin matrix $G$. Kernel functions like~\eqref{eq:fundsol} are
asymptotically smooth, i.e. there exist constants $C_{as} (\alpha,\beta)$ such that
\begin{equation}\label{eq:asymsmooth}
  \seminorm{\partial_x^\alpha \partial_y^\beta g(x,y)} \leq C_{as}(\alpha,\beta)
  \seminorm{x - y}^{- \seminorm{\alpha} - \seminorm{\beta}}
  \seminorm{g(x,y)},
  \quad x \neq y,
\end{equation}
for all multi-indices $\alpha, \beta \in \mathbb{N}^3$. Together with the
admissibility condition~\eqref{eq:qadm}, this property implies exponential
decay of the approximation error~\cite{MR2767920}.

\begin{theorem}[Approximation error]\label{thm:h2err}
  Let $r \times c \in \mathcal{P}^+$ be admissible with $\eta \in (0,2)$ and let
  $g(\cdot,\cdot)$ be an asymptotically smooth function. If we use a fixed number
  of $m$ interpolation points in each direction, resulting in $p = m^3$ points
  overall, the separable approximation
  \begin{equation*}
    \tilde{g}(x,y) = \sum_{\mu=1}^p \sum_{\nu=1}^p L_{r,\mu}(x)
    g(\xi_{r,\mu},\xi_{c,\nu}) L_{c,\nu}(y), \quad
    x \in X_r, y \in Y_c,
  \end{equation*}
  satisfies
  \begin{equation*}
    \norm{g - \tilde{g}}_{\infty,B_r \times B_c} \leq C_{int}(p)
    {\left(\frac{\eta}{2}\right)}^{m + 1} \norm{g}_{\infty,B_r \times B_c}
  \end{equation*}
  for some constants $C_{int}(p)$.
  Consequently, the matrix approximation error can be bounded by 
  \begin{equation*}
    \norm{G - \widetilde{G}}_F = \norm{
      \sum_{r \times c \in \mathcal{P}^+} G - U_r S_{r\times c} W_c^H
    }_F \leq C
    \max_{i \in I} \norm{\psi_i}_{L^2(\Gamma)}\,
    \max_{j \in J} \norm{\varphi_j}_{L^2(\Gamma)}\,
    {\left(\frac{\eta}{2}\right)}^{m + 1}
  \end{equation*}
  for some constant $C$ that depends on $C_{int}(p)$, $\Gamma$ and on the
  clustering.
\end{theorem}

\begin{algorithm}[H]
  \bgroup
  \everymath{\displaystyle}
  \caption{$\mathcal{H}^2$-matrix by interpolation.}\label{algo:h2matrix}
  \begin{algorithmic}[1]
    \Procedure{clusterbasis}{$r$, $p$}
    \If{$\sons(r) \neq \varnothing$}\Comment{Build cluster basis recursively}
    \For{$r^\prime \in \sons(r)$}
    \State{$E_{r^\prime,r}[\lambda,\mu] = L_{r,\mu}(\xi_{r^\prime,\lambda}),
      \quad i \in r^{\prime}, \mu = 1,\ \ldots,p$}\Comment{Transfer matrix}
    \State{${\left( U_{\hat{r}} \right)}_{\hat{r} \in \mathcal{L}(r^\prime)}$,
      ${\left( E_{r^*,\hat{r}} \right)}_
      {r^* \in \sons(\hat{r}),\, \hat{r} \in \mathcal{T}(r^\prime)}$
      = \Call{clusterbasis}{$r^\prime$,  $p$}}
    \EndFor%
    \Else\Comment{$r$ is leaf cluster, compute leaf matrix}
    \State{$U_r[i,\mu] = \int\limits_{\Gamma} \mathcal{L}_{r,\mu}(x) \psi_i (x)
      \mathop{dS(x)}, \quad i \in r,\ \mu = 1, \ldots,p$}
    \EndIf%
    \State{\textbf{return}
      ${\left( U_{r^\prime} \right)}_{r^\prime \in \mathcal{L}(r)}$,
      ${\left( E_{\hat{r},r^\prime} \right)}_{\hat{r} \in \sons(r^\prime),\,
        r^\prime \in \mathcal{T}(r)}$}
    \EndProcedure%
    \Statex%
    \Procedure{h2}{$b$}
    \If{$\sons(b) \neq \varnothing$}\Comment{Build $\mathcal{H}^2$-matrix
      recursively}
    \For{$b^\prime \in \sons(b)$}
    \State{$G[b^\prime]$ = \Call{h2}{$b^\prime$}}
    \EndFor%
    \Else%
    \If{$b = r \times c$ is admissible}\Comment{Compute coupling matrix}
    \State{$S_b[\mu,\nu] = g(\xi_{r,\mu}, \xi_{c,\nu}),
      \quad \mu,\nu = 1, \ldots, p$}
    \Else\Comment{Compute full matrix}
    \State{$G[i,j] = \int\limits_{\Gamma} \int\limits_{\Gamma} g(x,y)
      \varphi_j(y) \mathop{dS(y)} \psi_i(x) \mathop{dS(x)},
      \quad i \in r,\ j \in c$}
    \EndIf%
    \EndIf%
    \State{\textbf{return} $G[b]$}
    \EndProcedure%
  \end{algorithmic}
  \egroup
\end{algorithm}

\begin{figure}[hbt]
  \centering
  \hspace{0.589cm}\includegraphics[width=0.6\linewidth]{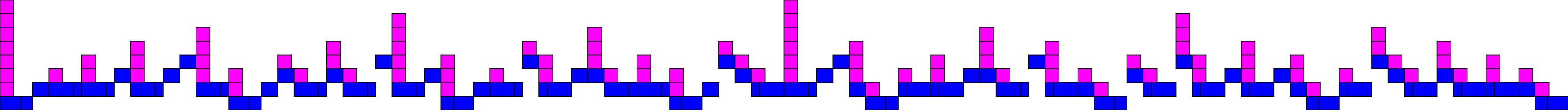}\\[1.53pt]
  \includegraphics[height=0.6\linewidth]{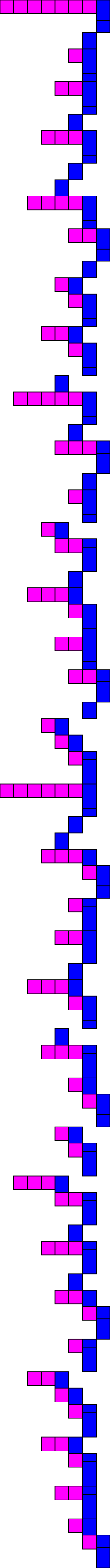}
  \includegraphics[width=0.6\linewidth]{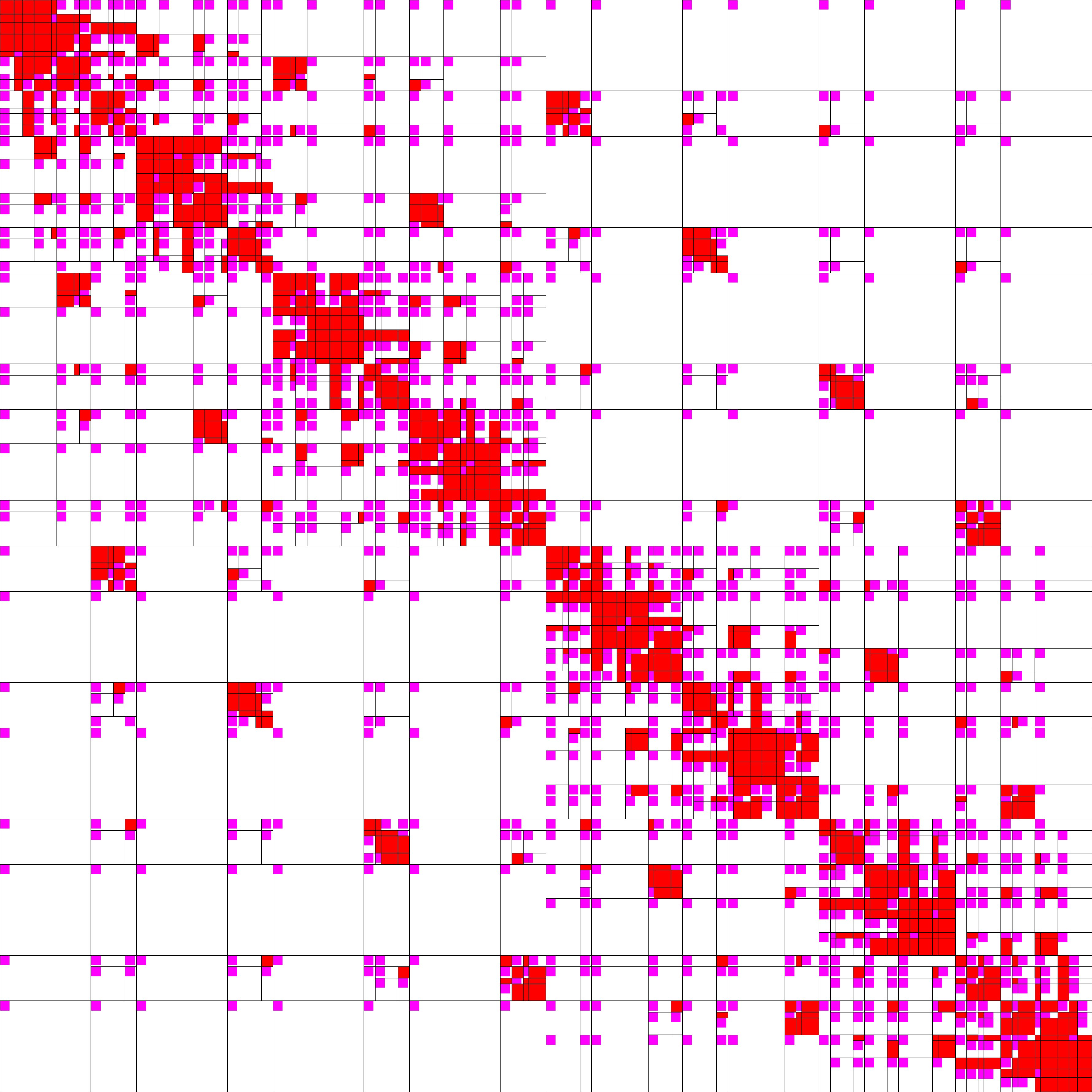}
  \caption{Visualisation of the storage costs of an
    $\mathcal{H}^2$-matrix.}\label{fig:h2matrix}
\end{figure}

As the computation of the far field only requires the assembly of the nested
cluster bases and coupling matrices, the storage costs are reduced drastically,
as depicted in Figure~\ref{fig:h2matrix}. The red boxes symbolise dense
near-field blocks, whereas far-field coupling matrices are painted magenta. The
blocks to the left and above the partition illustrate the nested row and column
cluster bases. There, leaf matrices are drawn in blue, while transfer matrices
are coloured in magenta.

The $\mathcal{H}^2$-matrix scheme scales linearly in the number of degrees of
freedom~\cite{MR2207953}.

\begin{theorem}[Complexity estimates]\label{thm:h2complex}
  Let $\mathcal{T}(I\times J)$ be sparse in the sense that a constant
  $C_{sp}$ exists such that
  \begin{equation*}
    \# \left\{ c^\prime \in \mathcal{T}(J) : r \times s^\prime
      \in \mathcal{T}(I\times J)\right\},
    \# \left\{ r^\prime \in \mathcal{T}(I) : r^\prime \times s
      \in \mathcal{T}(I\times J)\right\} \leq C_{sp}
  \end{equation*}
  for all $r \in \mathcal{T}(I)$ and $c \in \mathcal{T}(J)$.
  Then, the $\mathcal{H}^2$-matrix $\widetilde{G}$ requires
  \begin{equation*}
    \mathcal{O}(p(\# I + \#J))
  \end{equation*}
  units of storage and the matrix-vector multiplication can be performed in
  just as many operations.
\end{theorem}

The number of interpolation points $p=m^3$ equals the rank of the low-rank
factorisations in the far field and ultimately depends on the desired accuracy
$\varepsilon > 0$ of the approximation. We have
\begin{equation*}
  p = \mathcal{O}( {-\log (\varepsilon)}^{4})
\end{equation*}
in general, cf.~\cite{MR1993936}.

\begin{remark}
  Similar results hold for the kernel functions of the double layer and
  hyper-singular operator as well, see~\cite{MR2767920}.
\end{remark}

\section{Adaptive Cross Approximation}\label{sec:aca}
Returning to the setting of~\eqref{eq:tensor}, namely the approximation of the
tensor
\begin{equation*}
  \tensor{V}[i,j,k] = V_k[i,j],
\end{equation*}
we have the preliminary result that each slice $V_k$ is given in form of an
$\mathcal{H}^2$-matrix. Since the geometry $\Gamma$ is fixed for all
times, we can construct a partition that does not depend on the particular
frequency $s_k$. Therefore, we can select the same set of clusters
$\mathcal{T}(I)$ and $\mathcal{T}(J)$ for all $V_k$. In this way, the
partition $\mathcal{P}$ as well as the cluster bases ${\left(U_r \right)}_
{r \in \mathcal{T}(I)}$ and ${\left(W_c \right)}_{c \in \mathcal{T}(J)}$ have to
be built only once and are shared between all $V_k$. The latter only differ in
the coupling matrices and near-field entries, which have to be computed
separately for each frequency $s_k$.

Since all $V_k$ are partitioned identically, the tensor $\tensor{V}$ defined
in~\eqref{eq:tensor} inherits their block structure in the sense that it can be
decomposed according to $\mathcal{P}$ by simply ignoring the frequency index $k$.

\begin{definition}\label{def:tensopar}
  Let $K = \{0,\ldots,N-1\}$ and $\mathcal{T}(K) = \{ K \}$. In the current
  context, we define $\tensor{P} \in\mathcal{T}(I \times J \times K)$ to be the
  tensor partition with blocks
  \begin{equation*}
    b = r \times c \times K, \quad r \times c \in \mathcal{P},
  \end{equation*}
  which are admissible or inadmissible whenever $r \times c \in \mathcal{P}$ is
  admissible or inadmissible, respectively. 
\end{definition}
Naturally, this construction implies that the far-field blocks of $\tensor{V}$
are given in low-rank format,
\begin{equation*}
  \tensor{V}[r,c,k] = U_r\, S_{b,k}\, W_c^H, \quad k = 0, \ldots, N-1,
\end{equation*}
with $S_{b,k}$ being the coupling matrix of $b$ for the frequency $s_k$.
If we collect the matrices $S_{b,k}$ in the tensor $\tensor{S}_b$ in the same
manner as $V_k$ in $\tensor{V}$, we can factor out the cluster bases $U_r$ and
$W_c$ using the tensor product from Definition~\ref{def:tucker},
\begin{equation*}
  \tensor{V}[r,c,K] = \tensor{S}_b \times_1 U_r \times_2 \overline{W_c}.
\end{equation*}
The coupling tensor $\tensor{S}_b$ consists of kernel evaluations of the
transformed fundamental solution,
\begin{equation*}
  \tensor{S}_b[\mu,\nu,k] = \frac{e^{\textstyle{-s_k \seminorm{\xi_{c,\nu} - \xi_{r,\mu}}}}}
  {4 \pi \seminorm{\xi_{c,\nu} - \xi_{r,\mu}}},
\end{equation*}
which is smooth in $B_r \times B_c$ but also in the frequency
$s \in \mathbb{C}$. The latter holds true even for the near-field, whose entries
are
\begin{equation*}
  \tensor{V}[i,j,k] = \int\limits_\Gamma \int\limits_\Gamma
  \frac{e^{\textstyle{-s_k \seminorm{y - x}}}}
  {4 \pi \seminorm{y - x}} \varphi_j^0(y) \mathop{dS(y)} \varphi_i^0(x)
  \mathop{dS(x)}.
\end{equation*}
Hence, it is feasible to compress the tensor even further with respect to the
frequency index $k$. In particular, the above discussion shows that we may
proceed separately for each block $b \in \tensor{P}$, which represents either a
dense block $\tensor{V}[b]$ in the near-field or a coupling block $\tensor{S}_b$
in the far-field.

\subsection{Multivariate Adaptive Cross Approximation}\label{subsec:maca}
\begin{figure}[htb]
  \centering
  \includegraphics[width=0.7\textwidth]{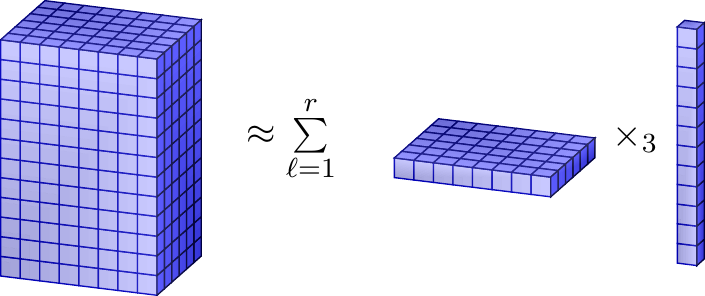}
  \caption{Visualisation of the low-rank factorisation.}\label{fig:lowrank_final}
\end{figure}

Let $\tensor{G} \in \mathbb{C}^{m \times n \times p}$ be a tensor. The
multivariate adaptive cross approximation (MACA) introduced in~\cite{MR2822767}
finds a low-rank approximation of rank $r \leq p$ of the form
\begin{equation}\label{eq:summaca}
  \tensor{G} \approx \tensor{G}^{(r)} = \sum_{\ell = 1}^r C_\ell \times_3 d_\ell
\end{equation}
with matrices $C_\ell \in \mathbb{C}^{m \times n}$ and vectors $d_\ell \in
\mathbb{C}^p$ as illustrated in Figure~\ref{fig:lowrank_final}. The main idea is
to reuse the original entries of the tensor.

\begin{figure}[htb]
  \centering
  \includegraphics[width=0.2\textwidth]{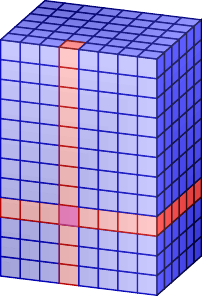}
  \includegraphics[width=0.431\textwidth]{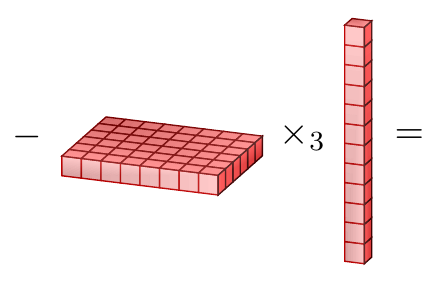}
  \includegraphics[width=0.2\textwidth]{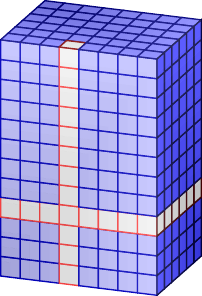}
  \caption{One step of the MACA.}\label{fig:maca_final}
\end{figure}

Starting from $\tensor{R}^{(0)} = \tensor{G}$, we pick a non-zero pivot element
in $\tensor{R}^{(\ell)}$ with index $(i_\ell, j_\ell, k_\ell)$ and select the
corresponding matrix slice and fibre for our next low-rank update, i.e.
\begin{equation*}
  C_\ell = \tensor{R}^{(\ell)}[1:m,1:n,k_\ell], \quad
  d_\ell = {\tensor{R}^{(\ell)}[i_\ell, j_\ell, k_\ell]}^{-1}
  \tensor{R}^{(\ell)}[i_\ell, j_\ell, 1:p].
\end{equation*}
Then, we compute the residual $\tensor{R}^{(\ell+1)}$ by subtracting their tensor
product,
\begin{equation*}
  \tensor{R}^{(\ell+1)} = \tensor{R}^{(\ell)} - C_\ell \times_3 d_\ell = 
  \tensor{R}^{(\ell)} -
  \frac{\tensor{R}^{(\ell)}[1:m,1:n,k_\ell] \times_3
    \tensor{R}^{(\ell)}[i_\ell, j_\ell, 1:p]}
  {\tensor{R}^{(\ell)}[i_\ell, j_\ell, k_\ell]}.
\end{equation*}
The residual $\tensor{R}^{(\ell+1)} = \tensor{G} - \tensor{G}^{(\ell)}$ measures
the accuracy of the approximation. After $r=\ell$ steps we obtain the low-rank
factorisation~\eqref{eq:summaca}. By construction, the cross entries successively
vanish, i.e.
\begin{equation*}
  \tensor{R}^{(r)}[i,j,k_\ell] = \tensor{R}^{(r)}[i_\ell,j_\ell,k] = 0, \quad
  \ell = 0, \ldots, r-1,
\end{equation*}
which implies $\tensor{R}^{(p+1)} = 0$ and hence $r \le p$.
Figure~\ref{fig:maca_final} depicts one complete step of the MACA. We extract
the cross consisting of $C_\ell$ and $d_\ell$ from $\tensor{R}^{(\ell)}$ and
subtract the update $C_\ell \times_3 d_\ell$, thereby eliminating the respective
cross from $\tensor{R}^{(\ell+1)}$.

The choice of the pivoting strategy is the crucial part of the algorithm. On the
one hand, it should lead to nearly optimal results, in the sense that high
accuracy is achieved with relatively low rank. On the other hand, it should be
reliable and fast, otherwise it would become a bottleneck of the algorithm.
Different pivoting strategies are available~\cite{MR2451321}, but we
restrict ourselves to finding the maximum entries in $C_\ell$ and $d_{\ell-1}$,
i.e. we choose $(i_\ell, j_\ell, k_\ell)$ such that
\begin{equation*}
  \begin{aligned}
    \seminorm{d_{\ell-1}[k_\ell]} &= \max_{k} \seminorm{d_{\ell-1}[k]},\\
    \seminorm{C_\ell[i_\ell,j_\ell]} &= \max_{i,j} \seminorm{C_\ell[i,j]},
  \end{aligned}
\end{equation*}
with $k_{1} = 0$. Throughout the algorithm, only $r$ slices and fibres of the
original tensor $\tensor{G}$ are used. Thus, there is no need to build the
whole tensor $\tensor{G}$ in order to approximate it and its entries are computed
only on demand. This feature presents a clear advantage of the ACA, especially
in BEM, where the generation of the entries is expensive. In this regard, the
routine \textproc{entry} in Algorithm~\ref{algo:maca} is understood to be a
call-back that computes the entries of $\tensor{G}$ at the time of its call.
Moreover, the tensors $\tensor{G}^{(\ell)}$ are never formed explicitly but are
stored in the low-rank format.

\begin{algorithm}[htb]
  \bgroup
  \everymath{\displaystyle}
  \caption{MACA}\label{algo:maca}
  \begin{algorithmic}[1]
    \Procedure{maca}{\textproc{entry}, $\varepsilon$}
    \State{$\tensor{G}^{(0)} = 0$, $k_{1} = 0$ and $\ell = 0$.}
    \Do%
    \State{$\ell = \ell + 1$}
    \State{$C_{\ell}[i,j] = \Call{entry}{i,j,k_{\ell}}-
      \tensor{G}^{(\ell-1)}[i,j,k_\ell], \quad i=1,\ldots,n,\ j=1,\ldots,m$}
    \State{$C_\ell[i_\ell,j_\ell] = \max_{i,j} \seminorm{C_\ell[i,j]}$}
    \If{$C_\ell[i_\ell,j_\ell] = 0$}
    \State{$\ell = \ell - 1$}
    \State{\textbf{break}}
    \EndIf%
    \State{$d_\ell[k] = {C_\ell[i_\ell,j_\ell]}^{-1} \left(
        \Call{entry}{i_\ell,j_\ell,k} - \tensor{G}^{(\ell-1)}[i_\ell,j_\ell,k]
        \right), \quad k = 1, \ldots, p$}
    \State{$\tensor{G}^{(\ell)} =
      \tensor{G}^{(\ell-1)} + C_\ell \times_3 b_\ell$}
    \State{$k_{\ell+1} = \argmax_{k}{\seminorm{d_{\ell}[k]}}$}
    \DoWhile{$\norm{C_\ell}_F \norm{d_\ell}_2 >
      \varepsilon \norm{\tensor{G}^{(\ell)}}_F$}
    \State{$r = \ell - 1$}
    \State{\textbf{return} $\tensor{G}^{(r)} =
      \sum_{\ell = 1}^r C_\ell \times_3 d_\ell$}
    \EndProcedure%
  \end{algorithmic}
  \egroup
\end{algorithm}

Here, we terminate the algorithm if the low-rank update
$C_\ell \times_3 d_\ell = \tensor{G}_\ell - \tensor{G}_{\ell-1}$ is sufficiently
small compared to $\tensor{G}^{(\ell)}$. Likewise, this stopping criterion does
not require the expansion of $\tensor{G}^{(\ell)}$ due to the identity
\begin{equation*}
  \begin{aligned}
    \norm{\tensor{G}^{(\ell)}}_F^2 &= \sum_{i,j,k}
    \seminorm{\sum_{\ell=1}^r C_\ell [i,j] b_\ell[k]}^2\\
    &= \sum_{\ell,\ell^\prime = 1}^r
    \left( \sum_{i,j} C_\ell[i,j]\, \overline{C_{\ell^\prime}[i,j]} \right)
    \left( \sum_{k} b_\ell[k]\, \overline{b_{\ell^\prime}[k]} \right).
  \end{aligned}
\end{equation*}
Neglecting the numerical work needed to compute the entries of $\tensor{G}$, the
overall complexity of the MACA amounts to $\mathcal{O}(r^2 (n m + p))$.

If we collect the vectors $d_\ell$ in the matrix $D^{(r)}\in
\mathbb{C}^{p\times r}$ and the matrices $C_\ell$ in the tensor $\tensor{C}^{(r)}
\in \mathbb{C}^{m \times n \times r}$, we obtain the short representation
\begin{equation}\label{eq:maca}
  \tensor{G}^{(r)} = \tensor{C}^{(r)} \times_3 D^{(r)},
\end{equation}
which is equivalent to~\eqref{eq:summaca}.

\begin{remark}\label{rem:unfolding}
  A tensor $\tensor{X} \in \mathbb{C}^{I_1 \times \cdots \times I_d}$ can be
  unfolded into a matrix by rearranging the index sets, which is called matricisation.
  For instance, the mode-$j$ unfolding $\mathcal{M}_j(\tensor{X}) \in
  \mathbb{C}^{I_j \times (\Pi_{k \neq j} I_k)}$ is defined by
  \begin{equation*}
    \mathcal{M}_j(\tensor{X})[i_j, (i_1,\ldots,i_{j-1},i_{j+1},\ldots,i_d)] =
    \tensor{X}[i_1,\ldots,i_d].
  \end{equation*}
  With this in mind, it turns out that the MACA is in fact the standard ACA
  applied to a matricisation of the tensor. In our special case, it is the mode-3
  unfolding.
\end{remark}

Due to Remark~\ref{rem:unfolding}, we can derive error bounds for the approximant
$\tensor{G}^{(r)}$ based on standard results for the ACA.
\begin{theorem}[Approximation error]\label{thm:acaerr}
  Let $\tensor{G}$ be either a dense block $\tensor{V}[b]$ or a coupling block
  $\tensor{S}_b$. Under the assumptions of~\cite[Theorem 3.35]{MR2451321}, there
  exist $0<\rho<1$ and $C>0$ such that the residual satisfies
  \begin{equation*}
    \norm{\tensor{R}^{(\ell)}}_F = \norm{\tensor{G} - \tensor{G}^{(\ell)}}_F
    < C \rho^{\ell+1},\quad \ell = 1,\ldots, N.
  \end{equation*}
  The constant $C>0$ depends on the block $b$ and on the distribution of the
  frequencies in the \textup{CQM}~\eqref{eq:cqm_freq}.
  \begin{proof}
    At first, we parameterise $s$ according to~\eqref{eq:cqm_freq},
    \begin{equation*}
      s(\theta) = \frac{\chi\left(R\cdot
        e^{\textstyle{\pi \imath \theta}}\right)}{\Delta t},
      \quad \theta \in [-1,1].
    \end{equation*}
    Then, the entries of $\tensor{G}$ are obtained by collocation of the
    functions
    \begin{equation*}
      \begin{aligned}
        F_{ij} (\theta) &= \int\limits_\Gamma \int\limits_\Gamma
        \frac{e^{\textstyle{{-s(\theta) \seminorm{y - x}}}}}
             {4 \pi \seminorm{y - x}} \varphi_j^0(y) \mathop{dS(y)}
        \varphi_i^0(x) \mathop{dS(x)}, \\
        G_{\mu \nu} (\theta) &=
        \frac{e^{\textstyle{-s(\theta) \seminorm{\xi_{c,\nu} - \xi_{r,\mu}}}}}
             {4\pi \seminorm{\xi_{c,\nu} - \xi_{r,\mu}}}
      \end{aligned}
    \end{equation*}
    at $\theta = k / N$. Since they are analytic in $\theta$,
    we may use~\cite[Section 68 (76)]{MR0048360} to bound the error of the
    best polynomial approximation of degree $\ell$,
    \begin{equation*}
      \inf_{m \in \mathcal{P}_\ell} \norm{f - m}_{\infty,[-1,1]} <
      \frac{2 M}{1 - \rho} \rho^{\ell+1},
    \end{equation*}
    where $f=F_{ij},G_{\mu \nu}$, $0<\rho<1$ and $M$ is chosen such that the
    absolute value of $f$ is less than $M$ within an ellipse in the complex plane
    whose foci are at $-1$ and $1$ and the sum of whose semi-axes is $1 / \rho$.
    Hence, the application of~\cite[Theorem 3.35]{MR2451321} yields the desired
    bound.
    \qed
  \end{proof}
\end{theorem}

Theorem~\ref{thm:acaerr} also justifies the choice of our stopping criterion.
If we assume
\begin{equation*}
  \norm{\tensor{R}^{(\ell+1)}}_F \le \delta \norm{\tensor{R}^{(\ell)}}_F,
\end{equation*}
then we obtain
\begin{equation*}
  \norm{\tensor{R}^{(r)}}_F \leq \delta \norm{\tensor{G}}_F
\end{equation*}
by setting $\varepsilon = \delta (1-\delta)/(1+\delta)$.

\section{Combined Algorithm}\label{sec:combalgo}
We are ready to state the complete algorithm, see Algorithm~\ref{algo:combi},
for the low-rank approximation of the boundary element tensors from
Section~\ref{subsec:galapx}. In the first step, we build the cluster bases and
construct a partition of the associated tensor~\eqref{eq:tensor} as outlined
in Section~\ref{sec:h2matrix} and Definition~\ref{def:tensopar}.
In the second step, we apply the MACA from Section~\ref{sec:aca} to each block of
the partition and obtain low-rank factorisations of the form~\eqref{eq:maca}.
Eventually, we end up with a hierarchical tensor approximation, which reads
\begin{equation}\label{eq:lrfinal}
  \begin{aligned}
    \tensor{V}[b] &\approx \tensor{C}_b \times_3 D_b, && b \in \tensor{P}^-,\\
    \tensor{V}[b] &\approx \tensor{C}_b \times_1 U_r \times_2 \overline{W_c}
    \times_3 D_b, && b \in \tensor{P}^+.
  \end{aligned}
\end{equation}
Besides the calls of the \textproc{MACA} routine, Algorithm~\ref{algo:combi} is
identical to Algorithm~\ref{algo:h2matrix}.

\begin{algorithm}[htb]
  \caption{Combined Algorithm}\label{algo:combi}
  \begin{algorithmic}[1]
  \Procedure{main}{${\left\{\psi_i\right\}}_{i \in I}$,
    ${\left\{\varphi_j\right\}}_{j \in J}$,
    ${\left\{s_n \right\}}_{n=1,\ldots,N}$,
    $n_{\min}$, $\eta$, $m$, $\varepsilon$}
    \State{$\mathcal{T}(I)$ = \Call{cluster}{${\left\{\psi_i\right\}}_{i \in I}$,
        $n_{\min}$}, \quad
      $\mathcal{T}(J)$ = \Call{cluster}{${\left\{\varphi_j\right\}}_{j \in J}$,
        $n_{\min}$}}
    \State{$rb$ = \Call{clusterbasis}{$I$, $m$},\quad
      $cb$ = \Call{clusterbasis}{$J$, $m$}}
    \State{$\mathcal{P}$ = \Call{partition}{$\mathcal{T}(I)$, $\mathcal{T}(J)$,
        $\eta$}}
    \For{$b \in \mathcal{P}$}\Comment{Call MACA for each block}
    \If{$b$ is admissible}
    \State{$\tensor{C}_b$, $B_b$ = \Call{maca}{\textproc{far}, $b$,
        $\varepsilon$}}
    \Else%
    \State{$\tensor{C}_b$, $B_b$ = \Call{maca}{\textproc{near}, $b$,
        $\varepsilon$}}
    \EndIf%
    \EndFor%
    \State{\textbf{return} $\tensor{A} = \left\{
      {\left\{\tensor{C}_b, B_b \right\}}_{b \in \mathcal{P}}, rb, cb \right\}$}
    \EndProcedure%
    \Statex%
    \Procedure{far}{$\mu$, $\nu$, $n$}
    \State{\textbf{return} $g(\xi_{r,\mu}, \xi_{c,\nu}, s_n)$}
    \Comment{Entries of coupling tensors}
    \EndProcedure%
    \Statex%
    \Procedure{near}{$i$, $j$, $n$}
    \State{\textbf{return} $\int\limits_\Gamma \int\limits_\Gamma g(x,y, s_n)
      \psi_j(y) \mathop{dS(y)} \varphi_i(x) \mathop{dS(x)}$}
    \Comment{Entries of dense blocks}
    \EndProcedure%
  \end{algorithmic}
\end{algorithm}

\subsection{Error analysis}\label{subsec:erroranalysis}
Before we discuss computational advantages of this algorithm, we briefly
state a basic result for the approximation error.

\begin{corollary}[Tensor Approximation Error]\label{cor:combierr}
  For every $\varepsilon \ge0$, we find an approximation $\widetilde{\tensor{V}}$
  of $\tensor{V}$ generated by Algorithm~\ref{algo:combi} which satisfies
  \begin{equation*}
    \norm{\tensor{V} - \widetilde{\tensor{V}}}_F  \le \varepsilon.
  \end{equation*}
  \begin{proof}
    For admissible blocks $b \in \tensor{P}^+$ we observe that the first term in
    \begin{equation*}
      \begin{aligned}
        \norm{\tensor{V}[b] - \tensor{C}_b \times_1 U_r \times_2 \overline{W_c}
          \times_3 D_b}_F &\le \norm{\tensor{V}[b] - \tensor{S}_b \times_1 U_r
          \times_2 \overline{W_c}}_F \\
        &\quad+\norm{\left(\tensor{C}_b \times_3 D_b - \tensor{S}_b \right)
          \times_1 U_r \times_2 \overline{W_c}}_F,
      \end{aligned}
    \end{equation*}
    is controlled by the $\mathcal{H}^2$-approximation and the second one by
    the MACA. By virtue of Theorems~\ref{thm:h2err} and~\ref{thm:acaerr}, we can
    prescribe accuracies $\delta_b >0$ on the approximation error for every
    far-field block. Similarly, we can bound the error block-wisely in the
    near-field by $\delta_b$. Hence, we obtain the desired bound by choosing
    $\delta_b$ such that
    \begin{equation*}
      \norm{\tensor{V} - \widetilde{\tensor{V}}}^2_F 
      \le \sum_{b \in \tensor{P}} \delta_b^2 \le \varepsilon^2
    \end{equation*}
    is satisfied.
    \hfill
    \qed
  \end{proof}
\end{corollary}

\begin{remark}
  The construction and analysis of low-rank approximations does not depend on
  the particular kernel function. Therefore, other boundary element matrices
  can be treated equivalently.
\end{remark}

\subsection{Complexity and Fast Arithmetics}\label{subsec:fastarith}

For admissible blocks, the low-rank approximation is given in the so called
Tucker format~\cite{MR3236394}.

\begin{definition}[Tucker format]\label{def:tucker}
  For a tensor $\tensor{X} \in \mathbb{C}^{I_1 \times \cdots \times I_d}$
  the Tucker format of tensor rank $(p_1,\ldots,p_d)$ consists of matrices
  $A^{(j)} \in \mathbb{C}^{p_j \times I_j}$, $j = 1, \ldots, d$, and a core
  tensor $\tensor{C} \in \mathbb{C}^{p_1 \times \cdots \times p_d}$ such that
  \begin{equation*}
    \tensor{X} = \tensor{C} \times_{j=1}^d A^{(j)},
  \end{equation*}
  with the tensor product from Definition~\ref{def:tenprod}. In the following,
  we call
  \begin{equation*}
    p = \max_{1\le j\le d} p_j
  \end{equation*}
  the maximum rank of the Tucker representation.
\end{definition}

One of the main advantages of the approximation in the Tucker format is the
reduction in storage costs. From~\eqref{eq:lrfinal}, we see that the Tucker
representation with maximum rank $p$ requires less than $\mathcal{O}(p^3 + p N +
p\,\#r + p\,\#c)$ units of storage compared with $\mathcal{O}(N\,\#r\,\#c)$ for
the dense block tensor. In addition, we improve the complexity even for
inadmissible blocks from $\mathcal{O}(N\, \#r\, \#c)$ to
$\mathcal{O}(p \,\#r\, \#c\, + p N)$. From these considerations, we immediately
deduce the following Corollary.

\begin{corollary}[Storage Complexity]\label{cor:complex}
  Under the assumptions of Theorem~\ref{thm:h2complex}, the hierarchical tensor
  decomposition needs about
  \begin{equation*}
    \mathcal{O}(p^2(\# I + \#J) + p N)
  \end{equation*}
  units of storage.
\end{corollary}

In addition, the low-rank structure allows us to substantially accelerate important
steps of the CQM. We recall that the computation of the integration weights
$\widehat{V}_n$ in~\eqref{eq:dft} comprises a matrix-valued discrete Fourier
transform of the auxiliary matrices $V_\ell$. If we use
representation~\eqref{eq:summaca} instead, we can factor out the
frequency-independent information, i.e.
\begin{equation}\label{eq:fastdft}
  \begin{aligned}
    \widehat{V}_n &= \frac{R^{-n}}{N} \sum_{\ell = 0}^{N-1}
    e^{\tfrac{-2\pi \imath}{N} n\ell} \sum_{k=1}^p C_k\, d_k[\ell] \\
    &=
    \sum_{k=1}^p C_k \frac{R^{-n}}{N} \sum_{\ell = 0}^{N-1}
    e^{\tfrac{-2\pi \imath}{N} n\ell} d_k[\ell], \quad n=0,\ldots,N.
  \end{aligned}
\end{equation}
Therefore, the transform has to be performed solely on the vectors $d_k$,
\begin{equation*}
  \hat{d}_k[n] = \frac{R^{-n}}{N} \sum_{\ell = 0}^{N-1}
  e^{\tfrac{-2\pi \imath}{N} n\ell} d_k[\ell],
\end{equation*}
with the result that the tensor of integration weights $\widehat{\tensor{V}}$
inherits the hierarchical low-rank format of the original tensor $\tensor{V}$.
In particular, the decomposition~\eqref{eq:lrfinal} still holds with $D_b$
replaced by $\widehat{D}_b$, whose columns are precisely the transformed vectors
$\hat{d}_k$. Thereby, we reduce the number of required FFTs from $\# r\, \# c$
to $p$ per block.
The other major improvement concerns the computation of the right-hand sides
in~\eqref{eq:cqmrhs}. There, discrete convolutions of the form
\begin{equation*}
  f_n = \sum_{k=0}^{n} \widehat{V}_{n-k}\, \underline{q}_k
\end{equation*}
need to be evaluated in each step. Once again, we insert~\eqref{eq:summaca}
and obtain
\begin{equation}\label{eq:fastrhs}
  f_n = \sum_{k=0}^{n} \sum_{\ell = 1}^r d_\ell[n-k] C_\ell\, \underline{q}_k =
  \sum_{\ell = 1}^p C_\ell \left(\sum_{k=0}^{n}d_\ell[n-k]\,
    \underline{q}_k \right).
\end{equation}
This representation requires in $\mathcal{O}(p)$ matrix-vector multiplications,
which amounts to $\mathcal{O}(p N)$ matrix-vector multiplications in total.
This is significantly less than the $\mathcal{O}(N^2)$ operations needed by
the conventional approach.

In combination with fast $\mathcal{H}^2$-matrix arithmetic~\cite{MR2207953}, the
algorithm scales nearly linearly in the number of degrees of freedom $M$ and time
steps $N$. This is shown Table~\ref{tab:complexity}, where we compare the storage
and operation counts of our fast algorithm with those of the traditional ones.
Note that the numerical effort for computing the tensor entries is not stated
explicitly but is reflected in the storage complexity. If efficient
preconditioners are available, we may replace the direct solver by an iterative
algorithm to eliminate the quadratic term $M^2$.

\begin{table}[htb]
  \centering
  \begin{tabular}{@{}ccccc@{}}
    \toprule
    & &\multicolumn{3}{c}{Computational} \\
    \cmidrule(r){3-5}
    Approximation & Storage & DFT & RHS & Solving \\
    \midrule
    None & $M^2 N$ & $M^2 N \log(N)$ & $M^2 N^2$ & $M^3 + M^2 N$\\
    $\mathcal{H}^2$ & $p M N$ & $p M N \log(N)$ & $p M N^2$ & $M^2 + M N$\\
    $\mathcal{H}^2$ + MACA & $p^2 M + p N$ & $p N \log(N)$ & $p^2 M N$
    & $M^2 + M N$ \\
    \bottomrule
  \end{tabular}
  \caption{Comparison of storage and computational complexity.
  }\label{tab:complexity}
\end{table}
  
\begin{remark}
  In~\cite{MR2452859}, an alternative method for solving~\eqref{eq:cqmgalapx} is
  presented. There, the convolutional structure is avoided by transforming to and
  from the Fourier domain. This approach leads to $N$ independent systems which
  involve the auxiliary matrices $V_k$ only. Thus, we may apply our
  approximation scheme also in this case to speed up the assembly of the matrices
  and reduce the number of matrix-vector multiplications needed to solve the
  systems.
\end{remark}

\section{Numerical Examples}\label{sec:numex}
In this section, we present numerical examples which confirm our theoretical 
results and show the efficiency of our new algorithm. In all experiments, we
set the parameter $\eta$ in the admissibility condition~\eqref{eq:admissibility}
to $2.0$ and choose $R = 10^{-5 / N}$ in the CQM.
The core implementation is based on the H2Lib software~\footnote{The source code
  is available at \url{https://github.com/H2Lib/H2Lib}.}.
The machine in use consists of two Intel Xenon Gold 6154 CPUs
operating at $3.00$~GHz with $376$~GB of RAM.

\subsection{Tensor Approximation}\label{subsec:tensor_approx}
The first set of examples concerns the performance and accuracy of the tensor
approximation scheme.

\begin{figure}[htb]
  \centering
  \includegraphics[width=0.6\linewidth]{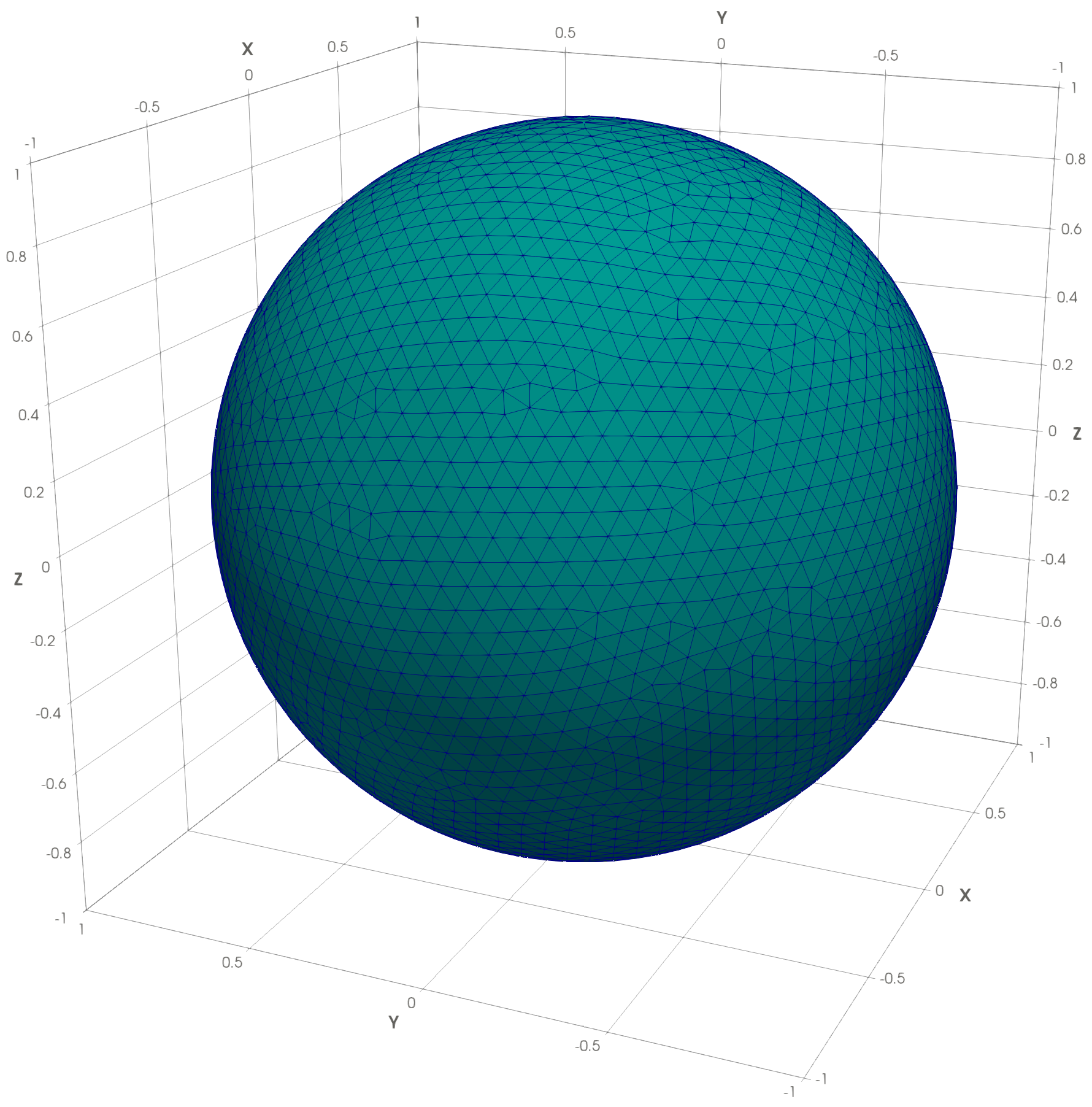}
  \caption{Spherical geometry $\Gamma$ used in the performance tests.
  }\label{fig:sphere}
\end{figure}

Let $\Gamma$ be the surface of a polyhedron $\Omega$, which approximates the
sphere of radius $1$ with $M$ flat triangles, see Figure~\ref{fig:sphere}. We
compare the dense tensor $\tensor{V}$ of single layer potentials with its
low-rank factorisation $\widetilde{\tensor{V}}$ and study the impact of the
interpolation order $m$ as well as accuracy $\varepsilon$ of the MACA on the
approximation error, rank distribution, memory requirements and computation time.

\begin{figure}[H]
  \centering
  \includegraphics{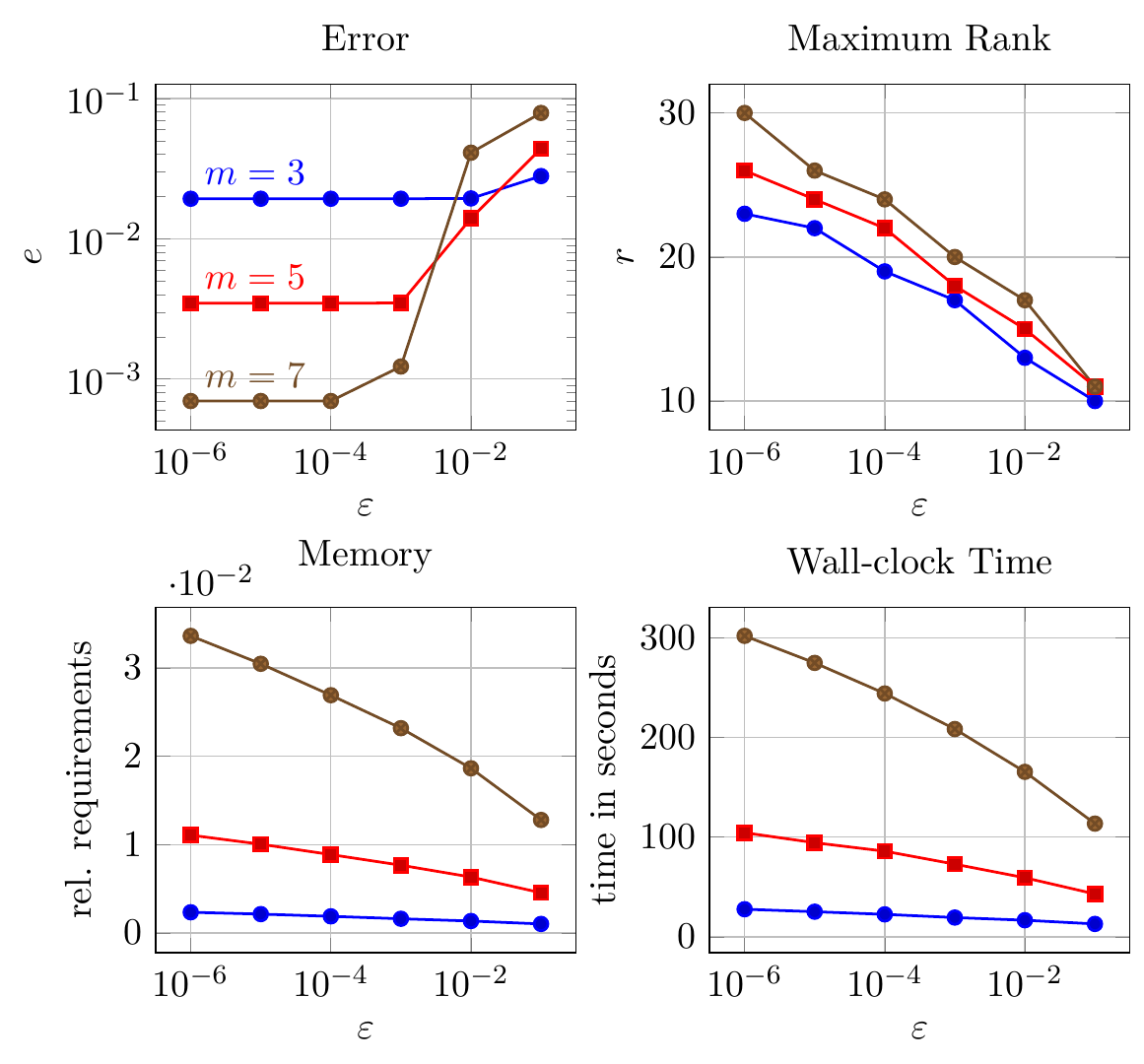}
  \caption{Results for $M=51200$, $N=256$, $m=3,5,7$ and
    improving accuracy $\varepsilon$ of the MACA.}\label{fig:eps}
  \includegraphics{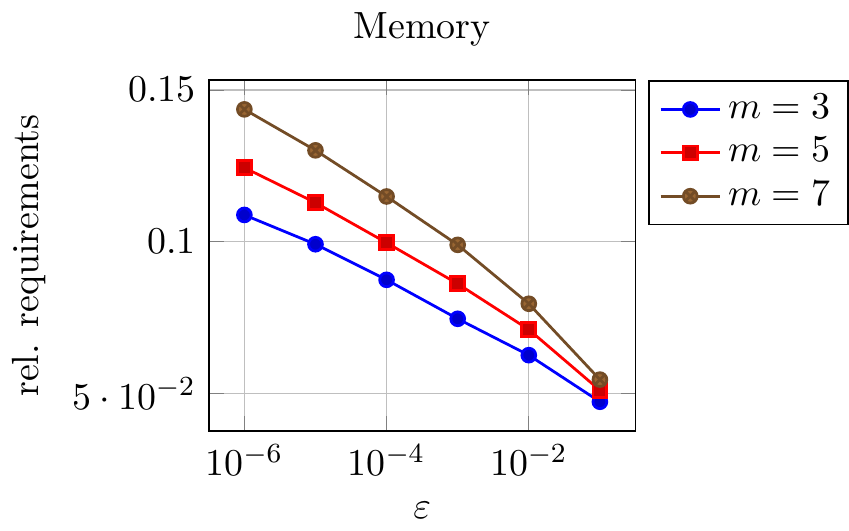}
  \caption{Storage requirements relative to the plain $\mathcal{H}^2$-matrix
    approach for $M=51200$, $N=256$, $m=3,5,7$ and varying $\varepsilon$.
  }\label{fig:memh2}
\end{figure}

We first set the number of degrees of freedom to $M = 51200$ and time steps
to $N=256$, resulting in a Courant number of $\Delta t / h \approx 0.83$.
In Figure~\ref{fig:eps}, the results for varying $\varepsilon$ and fixed
$m=3,5,7$ are presented. Foremost, we observe that the relative error
\begin{equation*}
  e = \frac{\norm{\tensor{V} - \widetilde{\tensor{V}}}_F}{\norm{\tensor{V}}_F}.
\end{equation*}
in the Frobenius norm decreases with $\varepsilon$ until it becomes constant
for high accuracies $\varepsilon \ge 10^{-4}$. This behaviour can be explained by
Corollary~\ref{cor:combierr}. Even if the coupling blocks are reproduced exactly
by the MACA, the $\mathcal{H}^2$-matrix approximation still dominates the total
error. Moreover, the numerical results confirm that the maximal block-wise rank
$r$ of the MACA depends logarithmically on $\varepsilon$ for fixed $m$. It stays
below $30$ in contrast to $256$ time steps, which reveals the distinct low-rank
character of the block tensors. Accordingly, our algorithm demands only for a
small fraction of memory compared with the conventional dense approach. At worst,
the compression rate reaches $3\%$ of the original storage requirements for
$m=7$. For $m=3,5$ and optimal choice of $\varepsilon = 10^{-2},10^{-3}$, we
further improve it to $0.2\%$ and $0.8\%$, respectively. Similarly, the
computation time needed for the assembly of the tensor is drastically reduced.
For the optimal values of $\varepsilon$, the algorithm takes only a couple
of seconds ($m=3,5$) or minutes ($m=7$) to compute the approximation of the
single layer potentials. Furthermore, we report that both memory requirements
and computation time scale logarithmically with $\varepsilon$.

To further demonstrate the benefits of the MACA, we consider the compression
rate in comparison to the case when only $\mathcal{H}^2$-matrices are used. In
Figure~\ref{fig:memh2}, the  storage costs for the same test setup are presented.
We notice that the inclusion of the MACA reduces the memory requirements to less
than $15\%$, while the same level of accuracy is achieved.

\begin{figure}[htb]
  \centering
  \includegraphics{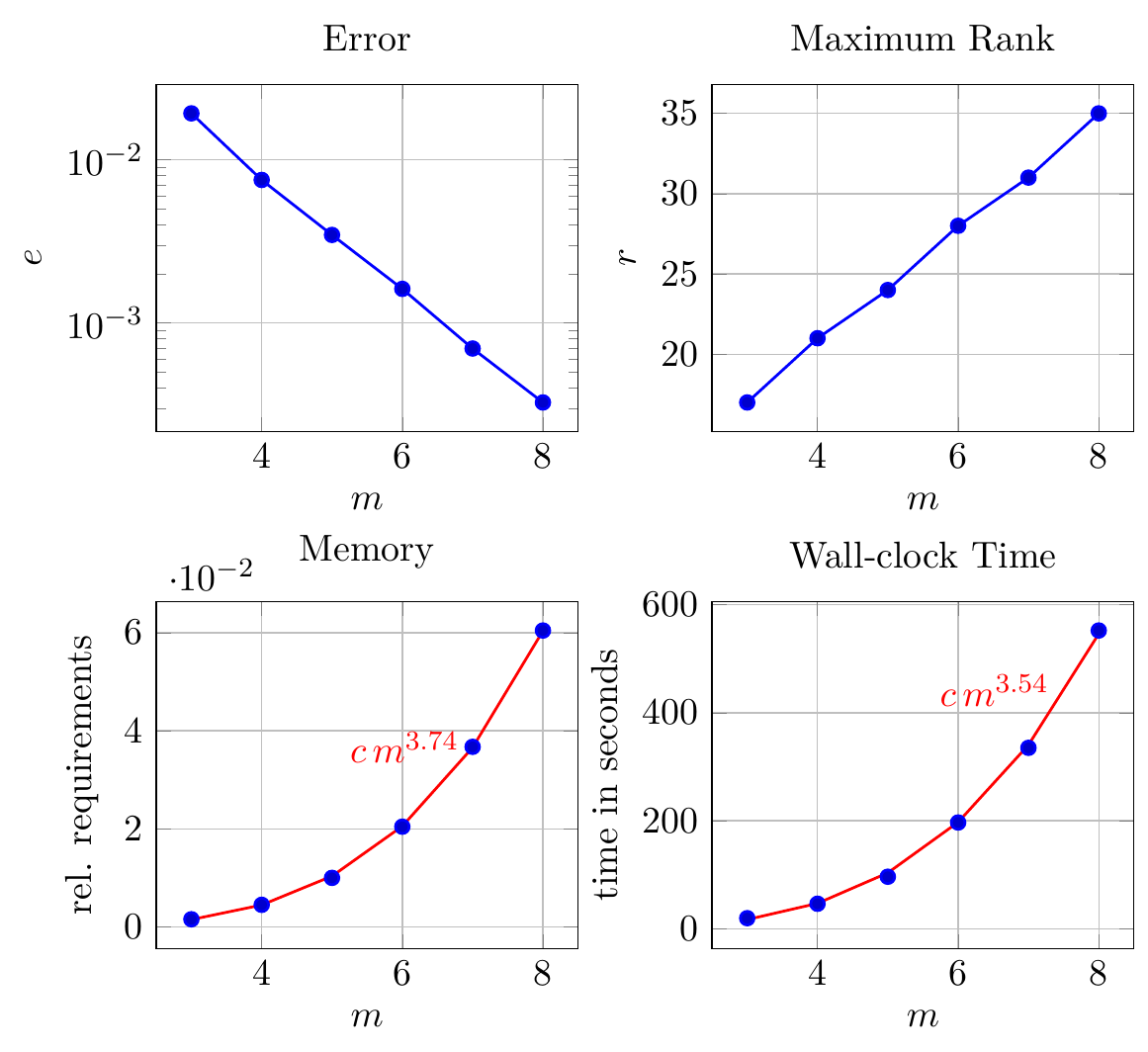}
  \caption{Results for $M=51200$, $N=256$, $\varepsilon=10^{-m}$ and increasing
    interpolation order $m$.}\label{fig:interp}
\end{figure}

If we modify the interpolation order $m$ instead, Theorems~\ref{thm:h2err}
and~\ref{thm:h2complex} indicate that the approximation error decreases
exponentially while the storage costs rise polynomially. This is confirmed by
the findings in Figure~\ref{fig:interp}, where we set $\varepsilon = 10^{-m}$ to
ensure that MACA error is negligible. Indeed, we see that the error $e$ is
roughly halved whenever $m$

\noindent
is increased by one and reaches almost $10^{-4}$ for
$m=8$. The upper right plot shows that the simultaneous change in $\varepsilon$
and $m$ leads to a linear growth of the MACA rank $r$ in terms of $m$. Since the
storage and computational complexity for fixed $M$ and $N$ is of order
$\mathcal{O}(p\, r)$, where $p = m^3$ is the interpolation rank, we observe
that the memory and time consumption scale approximately as $\mathcal{O}(m^4)$.
Note that in contrast to the prior example, the partition changes for every $m$,
since the latter directly affects the clustering.
We conclude that the algorithm yields accurate approximations with high
compression rates in short amounts of time. Moreover, we have seen that the
performance is more sensitive to a change in interpolation order $m$ than in
the MACA parameter $\varepsilon$. Hence, we recommend to select $\varepsilon$
on the basis of $m$ and not the other way around.

\begin{figure}[H]
  \centering
  \includegraphics{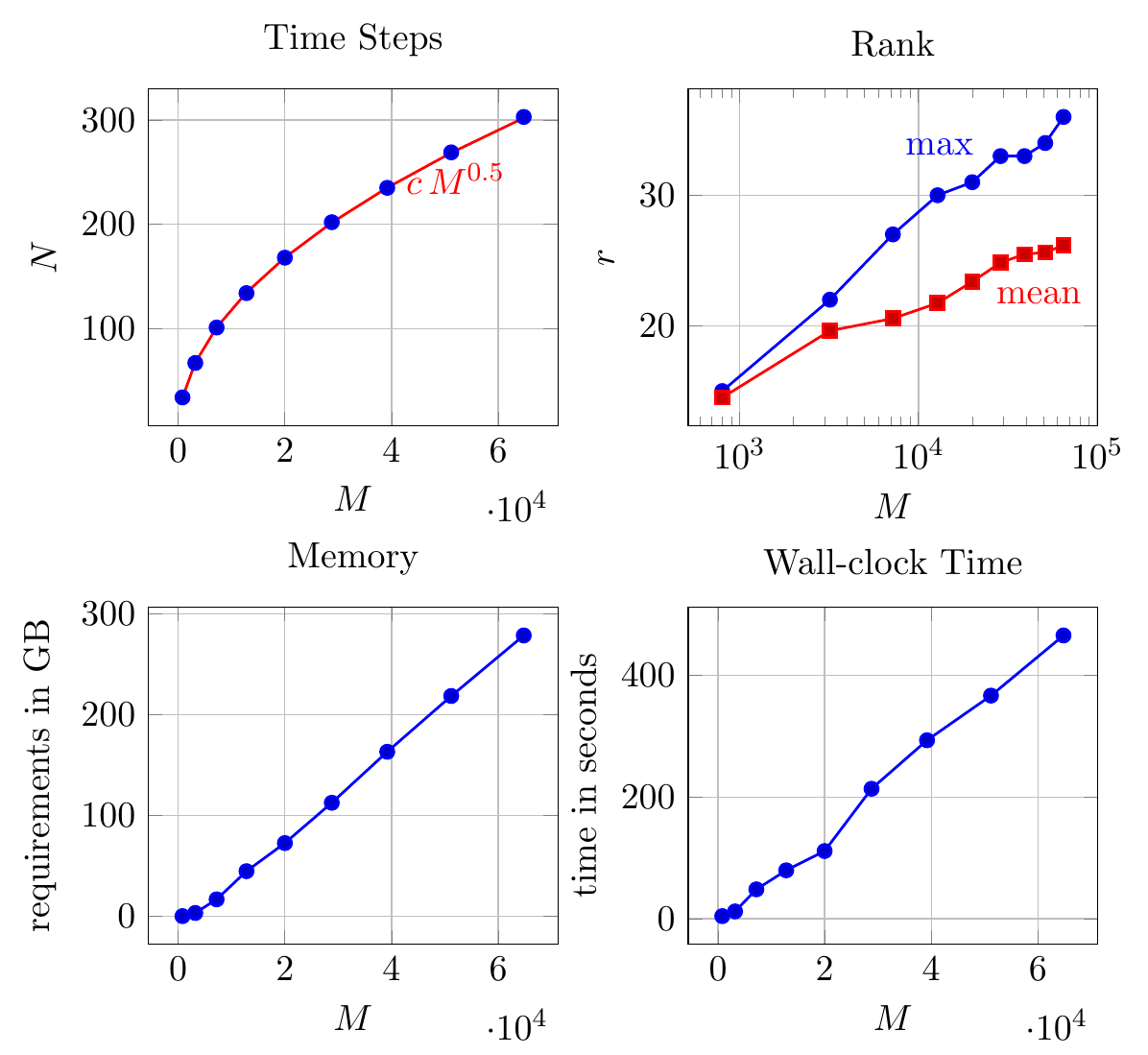}
  \caption{Results for constant courant number $\Delta t / h = 0.7$, $m=7$,
    $\varepsilon = 10^{-8}$ and increasing number of degrees of freedom $M$.
  }\label{fig:courant}
  \includegraphics[width=0.57\linewidth]{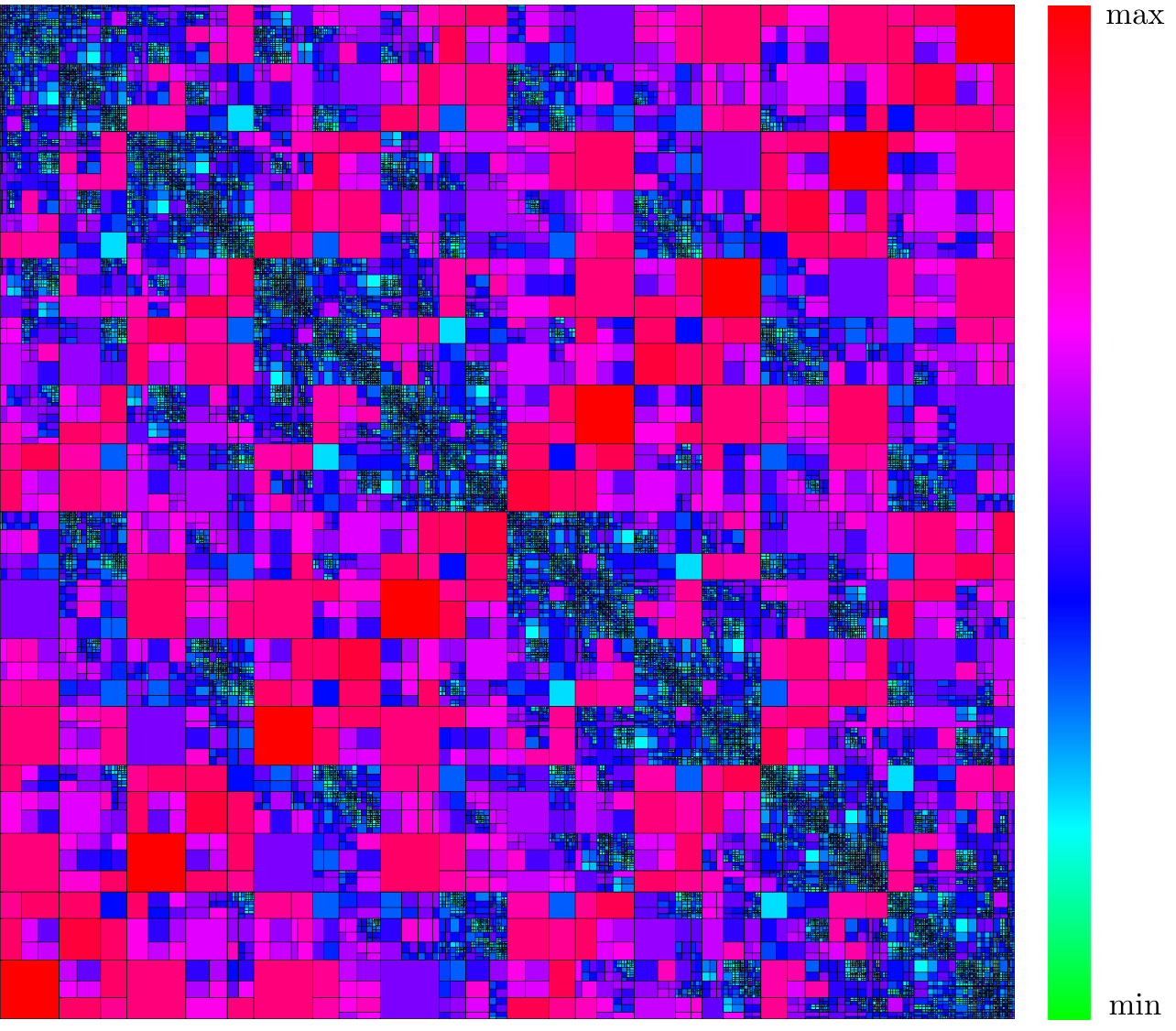}
  \caption{Exemplary rank distribution of the tensor partition for $M=51200$
    degrees of freedom and $N=768$ time steps.}\label{fig:rank}
\end{figure}

In the next two tests, we investigate the scaling of the algorithm in the number
of degrees of freedom $M$ and time steps $N$. First off, we fix the Courant
number $\Delta t / h = 0.7$ and refine the mesh $\Gamma$ successively. The
parameters $m=7$ and $\varepsilon = 10^{-8}$ are chosen in such a way that the
error $e$ is of the magnitude $10^{-3}$. Note that the approximation is more
accurate for small $M$ as the near-field still occupies a large part of the
partition. The results are depicted in Figure~\ref{fig:courant}, where the number
of time steps is added for the sake of completeness. First of all, we notice that
the storage and computational complexity are linear in $M$ in accordance with
Corollary~\ref{cor:complex}. Although the maximal rank $r$ grows logarithmically
at the same time, it does not influence the overall performance. This is probably
due to the average rank staying almost constant in comparison. The rank
distribution is visualised in form of a heat map in Figure~\ref{fig:rank}.

\begin{figure}[hbt]
  \centering
  \vspace{5.54976pt}
  \includegraphics{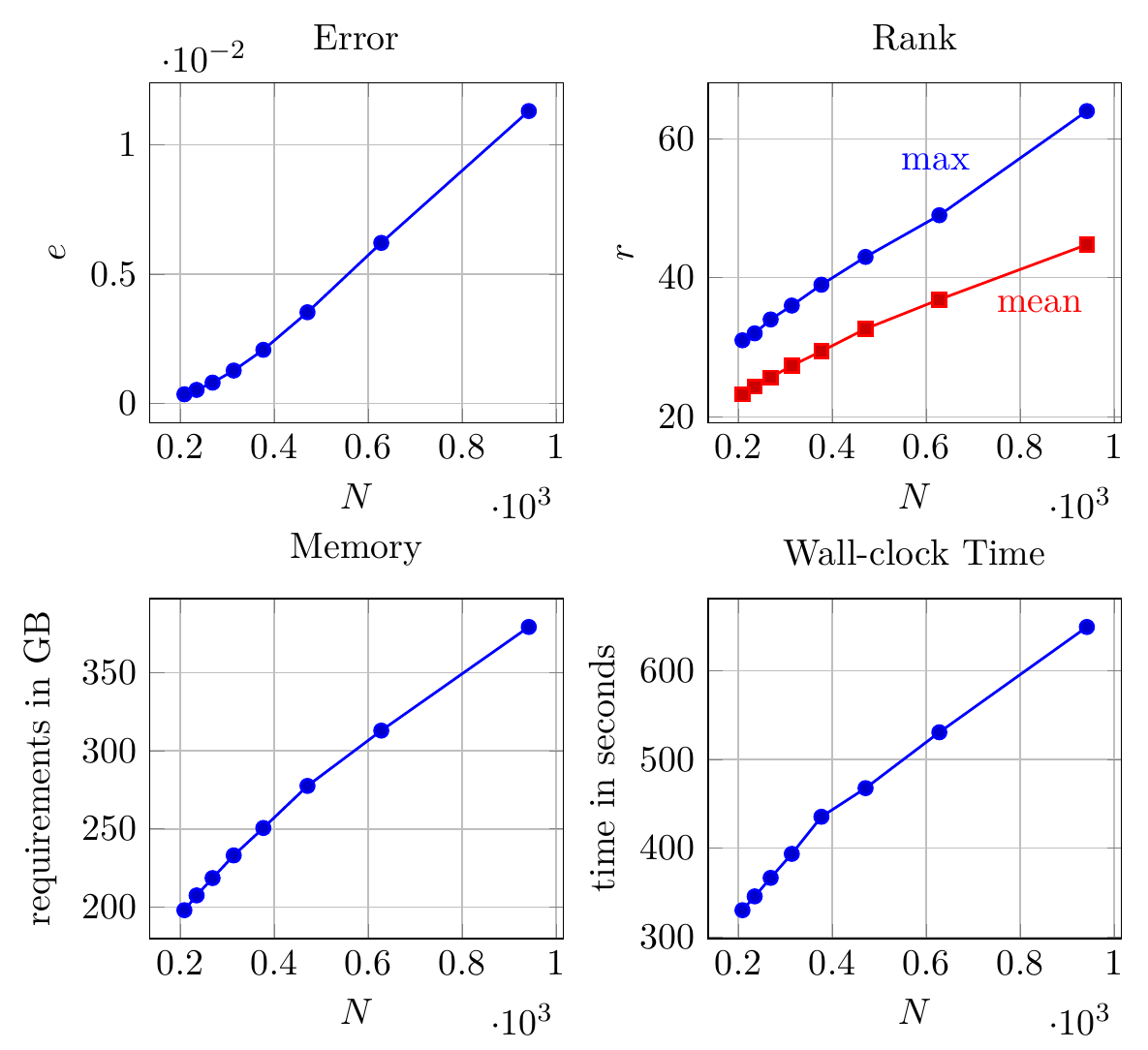}
  \caption{Results for  $M=51200$,  $m=7$, $\varepsilon = 10^{-8}$ and increasing
    number of time steps $N$.
  }\label{fig:timesteps}
\end{figure}

The dependence on the number of time steps $N$ for constant $M=51200$ is
illustrated in Figure~\ref{fig:timesteps}. We use the same values for the
parameters as before, i.e., $m=7$ and $\varepsilon = 10^{-8}$, and now change the
Courant number $\Delta t / h$ instead of $M$. The rise of the error $e$ is
attributed to the change in frequencies $s_\ell$. Since they grow in modulus, the
interpolation quality worsens with increasing $N$. We summarise that the
approximation scheme has linear complexity in both the number of degrees of
freedom $M$ and time steps $N$ for fixed tolerance $\varepsilon$ and
interpolation order $m$. 

\subsection{Scattering Problem}\label{subsec:scat_prob}
In this last section, we perform benchmarks for our fast CQM algorithm from
Section~\ref{sec:combalgo} and study the effect of the tensor approximation on
the solution of the wave problem.

\begin{figure}[htb]
  \centering
  \includegraphics[width=0.6\linewidth]{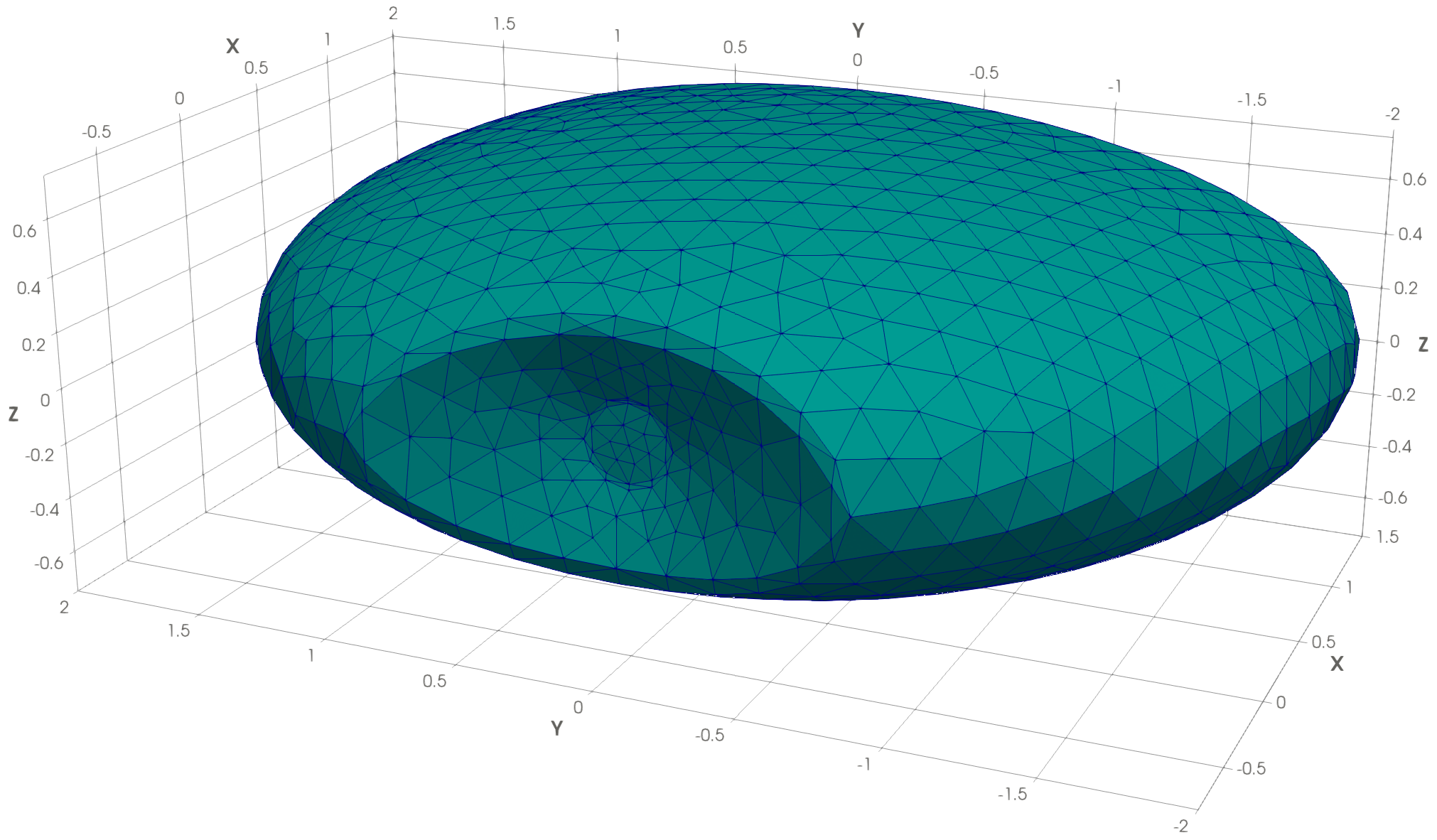}\\
  \includegraphics[width=0.6\linewidth]{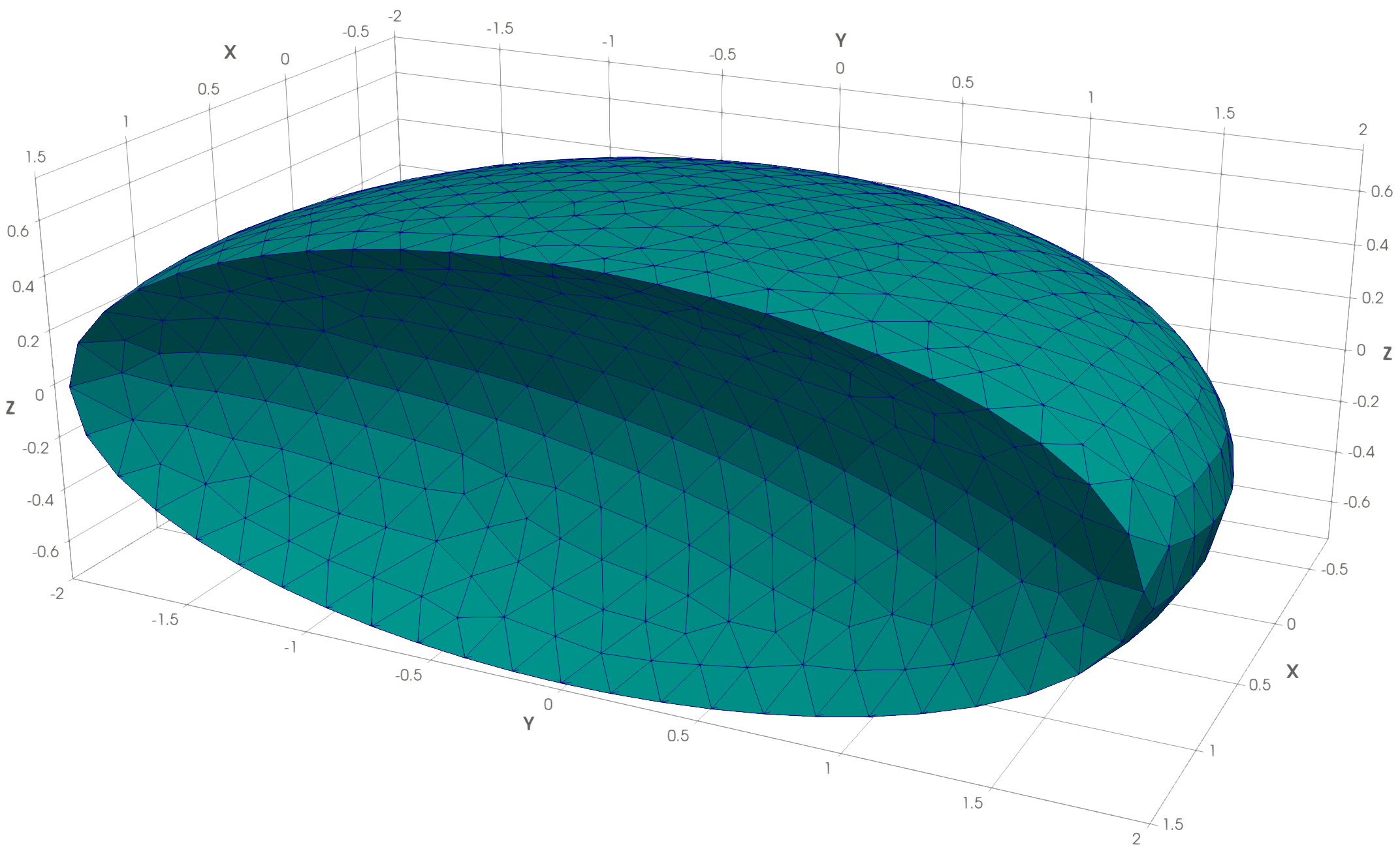}
  \caption{Front view and back view of the obstacle.}\label{fig:obstacle}
\end{figure}

To that end, we switch settings to the model problem~\eqref{eq:wave} posed in
the exterior of the geometry pictured in Figure~\ref{fig:obstacle}. The
spherical wave
\begin{equation*}
  u(x,t) = \frac{f(t - \seminorm{x})}{\seminorm{x}},
  \quad
  f(z) = \left\{
    \begin{aligned}
      &\cos(5 z + 1) - 1, && z > -1/5,\\
      &\,0, && z \leq -1/5,
    \end{aligned}
  \right.
\end{equation*}
serves as the exact solution. We shift the time variable such that $u$ reaches
the boundary $\Gamma$ right after $t=0$.

\begin{figure}[htb]
  \centering
  \includegraphics{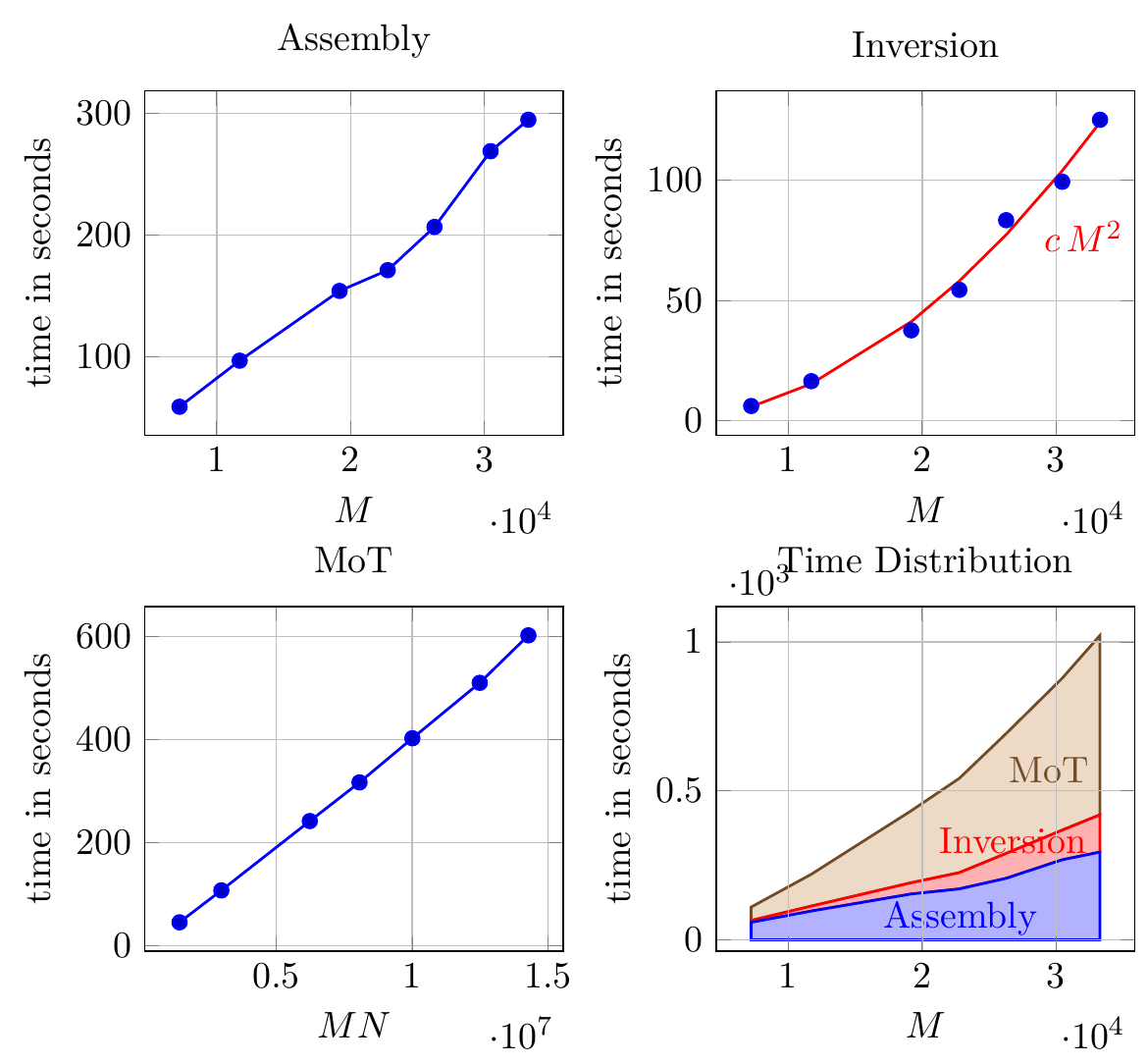}
  \caption{Results for constant courant number $\Delta t / h = 0.2$, $m=5$,
    $\varepsilon = 10^{-4}$ and increasing number of degrees of freedom $M$.
  }\label{fig:dir}
\end{figure}

The first part of tests concerns the fast arithmetics developed in
Section~\ref{subsec:fastarith}. To simplify matters, we impose pure Dirichlet
conditions on $\Gamma_D = \Gamma$ and solve for the Neumann trace $q_n$ in
\begin{equation}\label{eq:cqmdir}
  \widehat{V}_0 \underline{q}_n =
  \left(-\frac12 I + \widehat{K}_0\right) \underline{g}^D_n +
  \sum_{k=0}^{n-1} \left(\widehat{K}_{n-k} \underline{g}^D_k -
  \widehat{V}_{n-k} \underline{q}_k \right), \quad n=0,\ldots,N,
\end{equation}
as outlined in Section~\ref{subsec:galapx}. We choose $T=4.7$ as the final time.
We identify three major stages of the algorithm, firstly the assembly of the
tensors, secondly the inversion of $\widehat{V}_0$ and thirdly the step-by-step
solution of the linear systems, which is also known as marching-in-on-time (MoT).
In Figure~\ref{fig:dir}, we visualise how the running time is distributed among
the stages and how they scale in $N$ and $M$ for fixed parameters $m=5$ and
$\varepsilon=10^{-4}$. Overall, we see that the numerical results are consistent
with the estimates from Table~\ref{tab:complexity}. Beginning with the tensor
assembly, we once again observe linear complexity in $M$. The assembly comprises
the fast transformation from~\eqref{eq:fastdft} and is not explicitly listed,
since it requires less than $2$ seconds to perform in all cases. The LU
decomposition of the matrix $\widehat{V}_0$ involves $\mathcal{O}(M^2)$
operations but it nevertheless poses the least demanding part of the algorithm
for our problem size. On the other hand, the iterative solution
of~\eqref{eq:cqmdir} takes the largest amount of time. However, the application
of~\eqref{eq:fastrhs} allows for the fast computation of the right-hand sides in
just $\mathcal{O}(M N)$. This presents a significant speed up over the
conventional implementation. Taking into account that $N \sim M^{1/2}$ for a
constant Courant number, we expect the second stage to be the most expensive for
very large $M$. Therefore, it might be advantageous to switch to iterative
algorithms to solve the linear systems if efficient preconditioners are
available.

\begin{figure}[htb]
  \centering
  \includegraphics{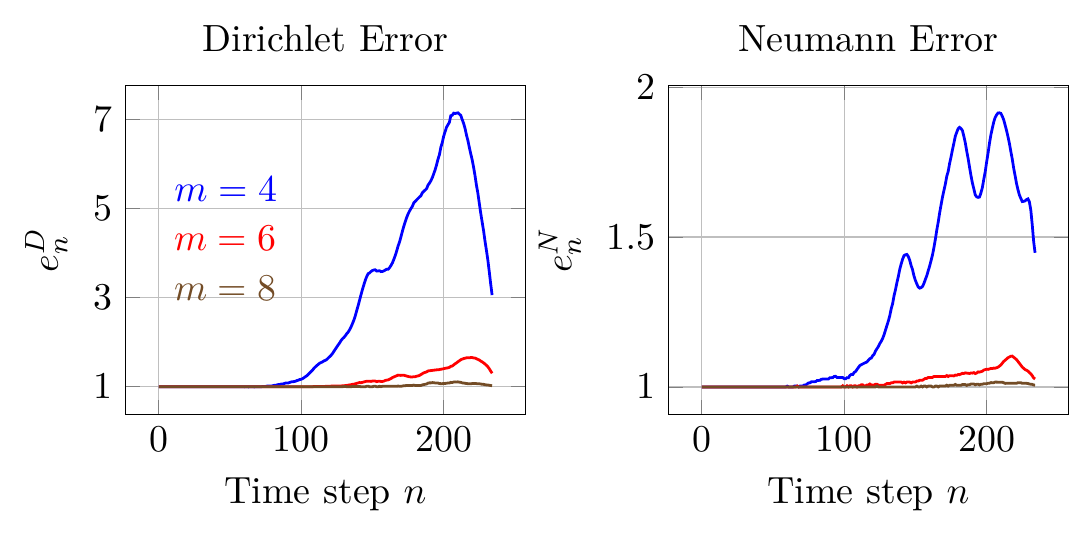}
  \caption{Results for constant $M=19182$, $N=235$, $\varepsilon = 10^{-8}$ and
    increasing interpolation order $m$.
  }\label{fig:mixed}
\end{figure}

For the last example, we split the boundary $\Gamma$ in Neumann and Dirichlet
parts $\Gamma_N$ and $\Gamma_D$ and replace the Dirichlet conditions by mixed
conditions. The Neumann boundary $\Gamma_N$ covers the upper half of $\Gamma$
with positive component $x_3>0$, while the rest of $\Gamma$ accounts to the
Dirichlet part. Furthermore, we consider the time interval $(0,1.7)$. We denote
by $u_n$ and $q_n$ the exact solutions and compare them with the approximations
$\tilde{u}_n$ and $\tilde{q}_n$ obtained by our fast CQM. We also include the
reference solutions $\tilde{u}^\textup{ref}_n$ and $\tilde{q}^\textup{ref}_n$
provided by the dense version. In particular, we are interested in how the
interpolation order $m$ affects the quality of the approximations,
which we estimate by computing the deviations
\begin{equation*}
  e^D_n = \frac{\norm{u_n - \tilde{u}_n}_{L^2(\Gamma)}}
  {\norm{u_n - \tilde{u}^\textup{ref}_n}_{L^2(\Gamma)}},
  \quad
  e^N_n = \frac{\norm{q_n - \tilde{q}_n}_{L^2(\Gamma)}}
  {\norm{q_n - \tilde{q}^\textup{ref}_n}_{L^2(\Gamma)}},
  \quad n=0,\ldots,N,
\end{equation*}
A value close to one indicates that the interpolation error does not spoil the
overall accuracy of the algorithm. We select a mesh with $M=19182$ triangles and
set the number of time steps to $N=235$. The choice of $\varepsilon=10^{-8}$
guarantees that the MACA does not deteriorate the interpolation quality. The
results for varying $m$ are depicted in Figure~\ref{fig:mixed}. We observe that
the approximations show the same level of accuracy for $n\le 75$ regardless of
the interpolation order. Then, the interpolation error becomes noticeable for
$m=4$ as $e^D_n$ and $e^N_n$ grow with $n$. The deviations from the reference
solution are considerably smaller for $m=6$ and our approximation scheme has
almost no impact on the accuracy for $m\ge 8$.

\section{Conclusion}\label{sec:conc}
In this paper, we have presented a novel fast approximation technique for the
numerical solution of wave problems by the CQM and BEM. We have given insights
into the theoretical and practical aspects of our algorithm and have explained
how it acts as a kind of hierarchical tensor approximation. Moreover, we have
proposed fast arithmetics for the evaluations of the discrete convolutions in
the CQM. In this way, we manage to reduce the complexity in terms of the number
of spatial degrees of freedom $M$ and number of time steps $N$: the storage costs
from $\mathcal{O}(M^2 N)$ to $\mathcal{O}(M + N)$ and the computational costs
from $\mathcal{O}(M^2 N^2 +M^3)$ to $\mathcal{O}(M N + M^2)$, where the quadratic
factor can be dropped if efficient preconditioners are available. Therefore, we
consider our work to be a step towards making large scale space-time simulations
possible with BEM.

\printbibliography%

\end{document}